\documentclass[Journal,letterpaper, InsideFigs, NoLineNumbers]{ascelike-new}
\usepackage[utf8]{inputenc}
\usepackage[T1]{fontenc}
\usepackage{lmodern}
\usepackage{graphicx}
\usepackage[figurename=Fig.,labelfont=bf,labelsep=period]{caption}
\usepackage{subcaption}
\usepackage{amsmath}
\usepackage{newtxtext,newtxmath}
\usepackage{algorithm}
\usepackage{algpseudocode}
\usepackage{multirow}
\usepackage[colorlinks=true,citecolor=blue,linkcolor=black]{hyperref}

\begin{document}

\title{Multi-Year Maintenance Planning for Large-Scale Infrastructure Systems: A Novel Network Deep Q-Learning Approach \thanks{This is a preprint of an article currently under review at the ASCE Journal of Infrastructure Systems.}} 

\author[1]{Amir Fard}
\author[1]{Arnold X.-X. Yuan}

\affil[1]{Department of Civil Engineering, Toronto Metropolitan University \\
350 Victoria Street, Toronto, ON, Canada, M5B 2K3}

\maketitle

\begin{abstract}
Infrastructure asset management is essential for sustaining the performance of public infrastructure such as road networks, bridges, and utility networks. Traditional maintenance and rehabilitation planning methods often face scalability and computational challenges, particularly for large-scale networks with thousands of assets under budget constraints. This paper presents a novel deep reinforcement learning (DRL) framework that optimizes asset management strategies for large infrastructure networks. By decomposing the network-level Markov Decision Process (MDP) into individual asset-level MDPs while using a unified neural network architecture, the proposed framework reduces computational complexity, improves learning efficiency, and enhances scalability.  The framework directly incorporates annual budget constraints through a budget allocation mechanism, ensuring maintenance plans are both optimal and cost-effective. Through a case study on a large-scale pavement network of 68,800 segments, the proposed DRL framework demonstrates significant improvements over traditional methods like Progressive Linear Programming and genetic algorithms, both in efficiency and network performance. This advancement contributes to infrastructure asset management and the broader application of reinforcement learning in complex, large-scale environments.
\end{abstract}

\section{Introduction}

Infrastructure asset management underpins the safety and reliability of critical systems that support everyday life, from road networks \cite{Chu2018,Yamany2024} and bridges \cite{Yang2022} to water distribution systems \cite{Rashedi2016} and sewer networks \cite{Fard2024}. Over time, these assets inevitably degrade and require maintenance, repair, or rehabilitation (MRR) measures. Public asset owners such as local and provicial governments or government agencies face the daunting task of selecting the best actions under restricted budgets \cite{Fard2024,Leppinen2025}, all while dealing with expanding infrastructure networks and the need for multi-year or even multi-decade planning \cite{Medury2014,Sasai2024}.

Effective maintenance strategies aim to optimize various performance indicators, including maximizing the Level of Service (LoS), minimizing life-cycle costs, and enhancing cost-effectiveness. The LOS usually is a network- or inventory-level performance measure characterizing how well the infrastructure assets as a whole have met users' service expectations and/or safety standards \cite{Fard2024}. Minimizing life-cycle costs involves reducing the total costs incurred over the assets' lifespan, encompassing both agency costs, such as expenses related to MRR interventions, and user costs associated with asset failure or poor performance \cite{Zhang2023}. Enhancing cost-effectiveness focuses on achieving the greatest improvement in asset conditions per unit of investment, ensuring that limited resources are utilized efficiently \cite{Yamany2024}.

However, current maintenance planning practices often adopt a myopic view, prioritizing immediate gains without considering long-term impacts. Such approaches may delay necessary maintenance, leading to accelerated deterioration and higher costs in the future. Therefore, there is a pressing need for methodologies that can optimize maintenance strategies over extended planning horizons, accounting for the dynamic nature of asset deterioration and the constraints imposed by budgets.

Budget constraints play a pivotal role in infrastructure asset management for they limit the available funding for maintenance activities in each planning period. When resources are limited, selecting a maintenance action for one asset affects the feasible actions for other assets due to the shared budget. This financial interdependency introduces combinatorial complexity, as the overall action space becomes the Cartesian product of individual asset actions constrained by the budget. Consequently, the optimization problem becomes more challenging, requiring methods that can efficiently handle large, interdependent decision spaces while satisfying budgetary limitations.

Historically, researchers have relied on methods such as linear programming \cite{Chu2018,Medury2014}, heuristic optimization (including evolutionary algorithms) \cite{Yamany2024,Elbeltagi2005,Fard2024}, and rule-based policies \cite{Rashedi2016} to devise maintenance plans. While effective for smaller-scale contexts or relatively simple problems, these techniques may not adapt well to networks containing tens of thousands of assets. They would also struggle to incorporate budget constraints in a way that balances immediate costs against long-term benefits and accounts for evolving asset conditions over extended planning horizons \cite{Fard2024,Sasai2024}.

Markov Decision Processes (MDPs) have gained attention as a robust theoretical framework for dealing with  the sequential nature of maintenance decisions \cite{Medury2014}. However, applying classical MDP approaches to large networks quickly becomes computationally infeasible \cite{Leppinen2025,Zhang2020} because the state and action spaces explode exponentially when each asset’s condition is considered. Storing and updating a value function for every possible combination of states and actions is simply impractical when asset numbers grow, leading to intense memory demands and computational overhead \cite{Rocchetta2019,Xu2024}.

Deep reinforcement learning (DRL) has recently emerged as a promising alternative, using deep neural networks to approximate value functions and policies \cite{Mnih2013}. Techniques such as Deep Q-Networks (DQN) allow an agent to learn which actions maximize future rewards through extensive experience, circumventing the need to keep track of explicit tables capturing every state-action pair \cite{Zhang2023}. DRL methods learn underlying patterns and can generalize to unseen states \cite{Zhou2022_2,Han2021}. Yet despite these breakthroughs, DRL-based infrastructure management solutions often face practical roadblocks: many studies are either demonstrated on relatively small or simplified networks, or they omit crucial constraints such as budgets and the interdependencies that arise from them \cite{Du2022,Rocchetta2019,Zhang2020,Han2021,Zhou2022_2}.

Moreover, some DRL approaches suffer from a significant increase in the dimensionality of their neural network input and output layers as the number of assets rises, potentially leading to intractable training times \cite{Zhang2020,Leppinen2025,Zhou2022}. Others provide only limited mechanisms for explicitly embedding budget or resource constraints into the learning process \cite{Du2022,Zhang2020,Han2021,Zhou2022_2}, leaving gaps in how real-world infrastructure agencies could implement these methods. For instance, it is not enough to identify which actions yield better future rewards if the algorithm does not also ensure that the recommended plan remains within strict financial limitations.

Meanwhile, multi-agent reinforcement learning (MARL) offers pathways to decompose large-scale problems by treating each component or asset as an individual agent. Value decomposition approaches like Value-Decomposition Networks (VDN) \cite{sunehag2017value} and QMIX \cite{rashid2020monotonic} have paved the way for centralized training with decentralized execution, making them appealing for cooperative settings. However, these methods mostly assume a shared reward without delving into budgetary constraints, and their factorization techniques are often based on specific additive or monotonic assumptions. More recent variants, such as QTRAN, Qatten, and FACMAC, seek to relax these constraints and improve coordination among agents \cite{son2019qtran,yang2020qatten,peng2021facmac}, but still do not fully tackle resource-limited scenarios on a large scale.

Researchers have begun exploring ways to incorporate resource constraints into MARL; for example, via constrained versions of QMIX (like CMIX) \cite{liu2021cmix} or using mean-field approaches for large agent populations \cite{mondal2024mean}. Others propose hierarchical methods that delegate constraint enforcement to specialized layers or decision-makers \cite{alqithami2025ch}. While all these efforts bring the field closer to a workable solution for complex, real-world problems, key gaps remain, particularly in scaling to thousands of assets, balancing local decisions against global resource limits, and capturing the intricate trade-offs within an integrated reinforcement learning framework.

Addressing these challenges motivates the development of a new DRL approach capable of bridging the gap between theoretical feasibility and practical deployment. This paper introduces a novel framework that decomposes a large infrastructure network into asset-level MDPs, each treated as an independent agent. By doing so, we mitigate the dimensionality explosion associated with tracking every asset concurrently at the network level. Importantly, we also embed a budget allocation component into the learning process, ensuring that maintenance actions across assets remain financially viable and cost-effective. In addition, by introducing a novel state representation scheme, our method employs a unified neural architecture in which the same model structure is shared among assets. This innovation not only enhances scalability but also leverages learning synergies across assets with similar characteristics.

The proposed state representation is comprehensive, blending asset-specific details (like current condition, estimated deterioration rate, and maintenance costs) with network-level features (such as overall LoS, the distribution of asset conditions, and historical budget usage). This hybrid perspective ensures that local decisions are informed by both local performance metrics and global targets. By systematically integrating budget considerations at every decision step, our framework can identify optimal strategies that effectively prioritize interventions to maintain a desired LoS or other performance metric while adhering to annual financial limitations.

In summary, our work develops a scalable DRL approach that reorganizes the large-scale maintenance problem into multiple subproblems and incorporates budgetary constraints directly into the learning and decision-making process. It employs a shared-parameter neural network, providing high learning efficiency across thousands of assets, and uses a state representation that integrates both local and global features to ensure alignment with network-wide objectives. We offer evidence of effectiveness through a case study with notable gains in LoS while annual budget constraints being honoured. By focusing on scalability and resource allocation, our framework not only enables infrastructure managers to prolong the lifespan and serviceability of critical networks under realistic budget scenarios, but also advances the practical application of reinforcement learning in complex, resource-constrained environments.

The remainder of this paper is structured as follows. In Section~\ref{sec:ProbFormulation}, we present the problem formulation, including the traditional optimization approach and the Markov Decision Process perspective, discussing different objective functions and the impact of budget constraints on problem complexity. Section~3 describes the proposed method, detailing the decomposition of the network MDP, the deep reinforcement learning framework, and the incorporation of budget constraints. In Section~4, we provide a detailed case study demonstrating the application of our method to a large-scale pavement network. Section~5 presents the results and discussion, highlighting the performance of our approach compared to traditional methods. Finally, in Section~6, we conclude with a discussion of the results and potential future research directions.

\section{Problem Formulation} \label{sec:ProbFormulation}

Maintenance planning for large-scale infrastructure systems involves selecting interventions to preserve or enhance asset performance over a multi-year horizon, all while complying with budget limitations. In this study, an infrastructure network is composed of \( n \) individual assets, each of which exhibits its own deterioration process, cost structure, and usage level or importance. Time is discretized into \( h \) discrete intervals (e.g., years), during which asset-specific maintenance decisions must be made. 

A maintenance plan specifies, for every asset and each planning interval, which treatment (if any) is to be applied. Possible treatments typically include a “No Intervention” option, along with one or more actions such as preventative maintenance, rehabilitation, or full reconstruction. Each action incurs a cost that depends on the asset’s size or importance measure and yields an expected improvement in condition. Although each asset’s deterioration model may be straightforward in isolation, combining many assets under a single multi-year budget constraint substantially complicates the decision problem.

To quantify how well the infrastructure network serves its users, many agencies use a network-level performance measure or Level of Service (LoS). The LoS aggregates individual assets’ conditions into a single representative indicator; for instance, a weighted sum of pavement quality indices for different roads. Maximizing the long-term LoS thus becomes the principal objective, although other objectives—such as minimizing total life-cycle costs or maximizing cost-effectiveness—may also be relevant. Regardless of the chosen objective, incorporating realistic budget constraints is essential because maintaining or improving one asset often leaves fewer resources available for other assets in the same funding cycle.

This section first outlines the classical mathematical formulation for multi-year maintenance planning, highlighting how the dimensionality grows with the number of assets and the need to enforce a shared annual budget. We then recast the problem as a Markov Decision Process (MDP), illustrating how decisions at each time step can be viewed through a sequential decision-making lens. The tension between local asset-level benefits and global resource constraints is made explicit, setting the stage for the methods proposed in subsequent sections.

\subsection{Mathematical Optimization Model}
A well-established way to formulate multi-year maintenance planning is through a nonlinear integer program that captures intervention decisions, system dynamics, and budget constraints:
\begin{align}
    \text{maximize} & \quad  \overline{\text{LoS}} = \frac{1}{h} \sum_{t=1}^{h} \text{LoS}_{t+1} = \frac{1}{h}  \sum_{t=1}^{h} \sum_{i=1}^{n} w_i \cdot \mathbb{E}[s_{i,t+1}], \label{eq:objective} \\
    \text{subject to} & \quad \nonumber \\
    & s_{i,t+1} = \sum_{j=0}^{m} x_{i,j,t} \cdot f_{i,j,t}(s_{i,t}), \quad \forall i = 1, \dots, n, \quad t = 1, \dots, h, \label{eq:state_transition} \\
    & \sum_{j=0}^{m} x_{i,j,t} = 1, \quad x_{i,j,t} \in \{0,1\}, \quad \forall i = 1, \dots, n, \quad t = 1, \dots, h, \label{eq:action_constraint} \\
    & \sum_{i=1}^{n} \sum_{j=1}^{m} c_{i,j,t} \cdot w_i \cdot x_{i,j,t} \leq b_t, \quad t = 1, \dots, h, \label{eq:budget_constraint}
\end{align}
where  \( x_{i,j,t} \) is a binary decision variable indicating whether treatment \( j \) is applied to asset \( i \) in year \( t \) (1 if applied, and 0 otherwise); \( b_t \) is the available budget in year \( t \); and \( f_{i,j,t}(s_{i,t}) \) is the transition function describing the effect of treatment \( j \) on asset \( i \) in year \( t \). When \( j = 0 \), it represents natural deterioration (no intervention).

Here the horizon-averaged Level of Service (HALoS) is used as the objective. When the end-of-horizon LoS (EHLoS) is used, the objective function is simply replaced by $\text{LoS}_h$. 

The transition functions \( f_{i,j,t}(s_{i,t}) \) are one-step prediction models that can be deterministic or stochastic. In a deterministic model, the future condition is a fixed function of the current condition and the applied action. This function can be a linear or nonlinear regression model based on asset attributes. In a stochastic model, the future condition is a random variable accounting for uncertainties in deterioration and maintenance effects. For example, a discrete-time Markov chain can be used, where the probability of transitioning from one condition state to another depends on the current state and the action applied.

Although all equations in the formulation all seems linear, one must note that the asset state transitions in Equation \eqref{eq:state_transition} are nonlinear functions of the binary decision variables in the preceding years through this recursive relationship.

As the number of assets \( n \), treatments \( m \), and the planning horizon \( h \) increase, the problem's complexity grows exponentially. The total number of possible action combinations is \( (m+1)^{nh} \) (or \( 2^{(m+1)nh} \) in its binary form), making it computationally infeasible to exhaustively search for the optimal solution in large-scale networks. This combinatorial explosion significantly challenges optimization efforts, particularly as the number of variables and constraints increases.

To address these challenges, population-based and heuristic algorithms, such as genetic algorithms, have been widely used in the literature. However, \citeN{Fard2024} highlights that these heuristic methods struggle to handle the complexity and constraints inherent in large-scale infrastructure management problems. While these approaches can yield satisfactory solutions for smaller networks, they often produce suboptimal or infeasible solutions when applied to real-world, large-scale networks. The nonlinear relationships among decision variables and the size of the solution space make it challenging for these algorithms to converge to optimal solutions within a reasonable time frame. Consequently, new methods, including reinforcement learning frameworks, are being explored to overcome these limitations and provide more effective solutions for large-scale asset management.

\subsection{Markov Decision Process Perspective}
\label{subsec:mdp_perspective}

The multi-year maintenance planning problem described above can be viewed through the lens of a Markov Decision Process (MDP), as the MDP provides a natural modelling framework for sequential decision-making under uncertainty \cite{Sutton18}. Within this framework, each year corresponds to a decision epoch where maintenance actions (interventions) are selected, and the condition states of assets evolve according to known or estimated deterioration dynamics. This subsection first describes an asset-level MDP formulation and then extends the idea to a network-level MDP under budget constraints, highlighting the ensuing computational challenges.

\subsubsection{Asset-Level MDPs}

Consider a single infrastructure asset \(i\). Let \(\mathcal{S}_i\) be the set of its possible condition states (e.g., condition index values), and let \(\mathcal{A}_i=\{0,1,\dots,m\}\) be the set of actions, where \(0\) represents ``No Intervention'' and \(1,\dots,m\) denote various maintenance treatments (e.g.\ rehabilitation, reconstruction). An asset-level MDP can be defined by $\bigl(\mathcal{S}_i,\;\mathcal{A}_i,\;P_i,\;R_i,\;\gamma\bigr)$, where \(\gamma\in(0,1]\) is a discount factor. At decision epoch \(t\), the asset is in state \(s_{i,t}\in\mathcal{S}_i\). An action \(a_{i,t}\in\mathcal{A}_i\) is applied, leading to a new state \(s_{i,t+1}\) according to the transition function  $
P_i\bigl(s_{i,t+1}\,\big\vert\,s_{i,t},\,a_{i,t}\bigr) = f_{i,a_{i,t},t}(s_{i,t})$.

The immediate reward $ r_{i,t} \;=\; R_i\bigl(s_{i,t},\,a_{i,t},\,s_{i,t+1}\bigr) $ captures the benefit (or negative cost) of transitioning from \(s_{i,t}\) to \(s_{i,t+1}\) under action \(a_{i,t}\). Design of the reward function often depends on the underlying objective of the original optimization problem. For instance, if the goal is to minimize life-cycle cost (LCC), one might define \(\,r_{i,t}\equiv -\,(\text{intervention cost}+\text{user cost})\).  In a simpler unconstrained objective of maximizing the asset’s state, one might set  $ R_i(s_{i,t},\,a_{i,t},\,s_{i,t+1}) \;=\; s_{i,t+1}$  assuming a higher condition state yields a higher reward. As it will be shown later, however, budget constraints often complicate this approach, especially when budget resources must be shared across multiple assets.  Nevertheless, with the immediate reward defined, the backward cumulative discounted return from time $t$ to the end of the planning horizon $h$ is expressed as 
\[
G_t
\;=\;
\sum_{\tau=t}^{h}
\gamma^{\tau-t}\,
r_{i,\tau} .
\]
When $r_{i,\tau} = s_{i,t}$ and $\gamma =1$, the expected cumulative return at time 0 equals the original objective function defined in Eq.~(\ref{eq:objective}). 

An MDP aims to find an optimal \emph{policy} that maximizes the expected cumulative discounted reward. A policy $\pi$ is a mapping from the state set $\mathcal{S}_i$ to the action set $\mathcal{A}_i$, and can generally be expressed as a policy function denoted by \(\pi_i(a\mid s)\), representing the probability of selecting action \(a\) in state~\(s\).  Under a given policy $\pi$ and a given state at time $t$, a \emph{value function} can be defined as
\[
V_i^\pi(s)
\;=\;
\mathbb{E}_\pi\!\Bigl[
\,\sum_{\tau=t}^{h} 
\gamma^{\,\tau - t}\,r_{i,\tau}
\;\Big\vert\;s_{i,t}=s
\Bigr] = \;
\mathbb{E}_\pi\!\Bigl[
G_t
\;\Big\vert\;s_{i,t}=s
\Bigr].
\]
With a further condition that the action at time $t$ is chosen as $a$, then an \emph{action-value function} can be defined as
\[
Q_i^\pi(s,a)
\;=\;
\mathbb{E}_\pi\!\Bigl[
G_t
\;\Big\vert\;s_{i,t}=s,\;a_{i,t}=a
\Bigr].
\]
The optimal value function \(V_i^*\) and optimal action-value \(Q_i^*\) satisfy Bellman optimality equations \cite{Sutton18}. At the asset level, an optimal policy \(\pi_i^*\) greedily selects
\[
\pi_i^*\bigl(a\mid s\bigr)
\;=\;
\arg\max_{a\in\mathcal{A}_i}\,
Q_i^*(s,a).
\]

\subsubsection{Network-Level MDP Under No Constraints}
For a network maintenance optimization without budget constraints, the asset-level MDP can be readily expanded to obtain a network-level MDP. In this case, the \emph{network state} at time \(t\) is simply augmented to a vector, i.e., $\mathbf{s}_t \;=\; \bigl(s_{1,t},\dots,s_{n,t}\bigr)$, while a \emph{joint action}  $\mathbf{a}_t
=\bigl(a_{1,t},\dots,a_{n,t}\bigr) $ is similarly defined. If each asset’s state evolves independently when its action is chosen, then the state transition of the network state is characterized by the product of the asset-level probability transition relations, i.e.,
\begin{equation}
\label{eq:network_transition_independent}
P\bigl(\mathbf{s}_{t+1}\,\big\vert\,\mathbf{s}_t,\mathbf{a}_t\bigr)
\;=\;
\prod_{i=1}^n 
P_i\bigl(s_{i,t+1}\,\big\vert\,s_{i,t},\,a_{i,t}\bigr).
\end{equation}
Similarly, when network effects are neglected, the network-level immediate reward function is a simple or weighted sum of the asset-level rewards: 
\[
R\bigl(\mathbf{s}_t,\mathbf{a}_t\bigr)
\;=\;
\sum_{i=1}^n
w_i \, R_i\bigl(s_{i,t},\,a_{i,t},\,s_{i,t+1}\bigr) .
\]
Without \emph{unconstraints} that would confound the feasibility of the asset-level optimal actions, one can readily show that the network-level MDP decomposes into \(n\) independent asset-level MDPs:
\[
Q^*(\mathbf{s},\,\mathbf{a})
\;=\;
\sum_{i=1}^n
Q_i^*\!\bigl(s_i,\,a_i\bigr),
\]
and each asset’s optimal action can be selected separately. This drastically reduces complexity, avoiding the exponential blow-up from joint actions.

\subsubsection{Budget Constraints and the Constrained MDP Challenge}

However, real-world maintenance planning imposes a shared annual budget \(b_t\) that couples each asset’s decision in year~\(t\), because $\mathbf{a}_t$ is feasible only if $ \sum_{i=1}^n
\bigl(w_i\,c_{i,a_{i,t},t}\bigr)
\,\le\,
b_t$. In other words, this single constraint disrupts the neat decomposition in the unconstrained network MDPs, because  if asset~\(i\) chooses an expensive intervention, fewer resources can be allocated to other assets. Therefore, although one can still formally write a network-level MDP with state \(\mathbf{s}_t\) and action \(\mathbf{a}_t\),  the budget constraint multiplies the dimensionality of feasible joint actions and invalidates the additive reward decomposition. 

Solving the constrained MDP exactly is typically intractable for large \(n\), since there may be \((m+1)^n\) joint actions at each step, and an exact dynamic program would require enumerating all feasible subsets of assets each year. Additionally, even approximate dynamic programming techniques must contend with the combined burden of high dimensional state \(\mathbf{s}_t\) and a combinatorial feasibility condition on \(\mathbf{a}_t\). This exponential growth is often referred to as the \emph{curse of dimensionality in the constrained setting}.

To recap, a network MDP without budget constraints can be completely decomposed into independent asset-level MDPs; local solutions at the asset level can be simply aggregated into a global optimum. Once budget constraints enter the picture, a single asset’s intervention can block or reduce feasible interventions for others. As a result,  the asset-level optimal policies often overestimate the returns from high-cost actions if the budget is insufficient to implement them for all assets simultaneously. Realistic approaches must explicitly consider the budget constraints either by introducing resource-allocation subroutines or by designing specialized multi-asset algorithms. In the next section, we propose a framework that (i) maintains the local MDP structures for each asset, (ii) uses a budget-allocation layer to select feasible actions under annual constraints, and (iii) learns policies that blend cost-normalized action-values at the asset level with global performance objectives in a scalable manner.

\section{Proposed Method}
\label{sec:proposed_method}

\begin{figure}[H]
    \centering
    \includegraphics[width=0.92\textwidth]{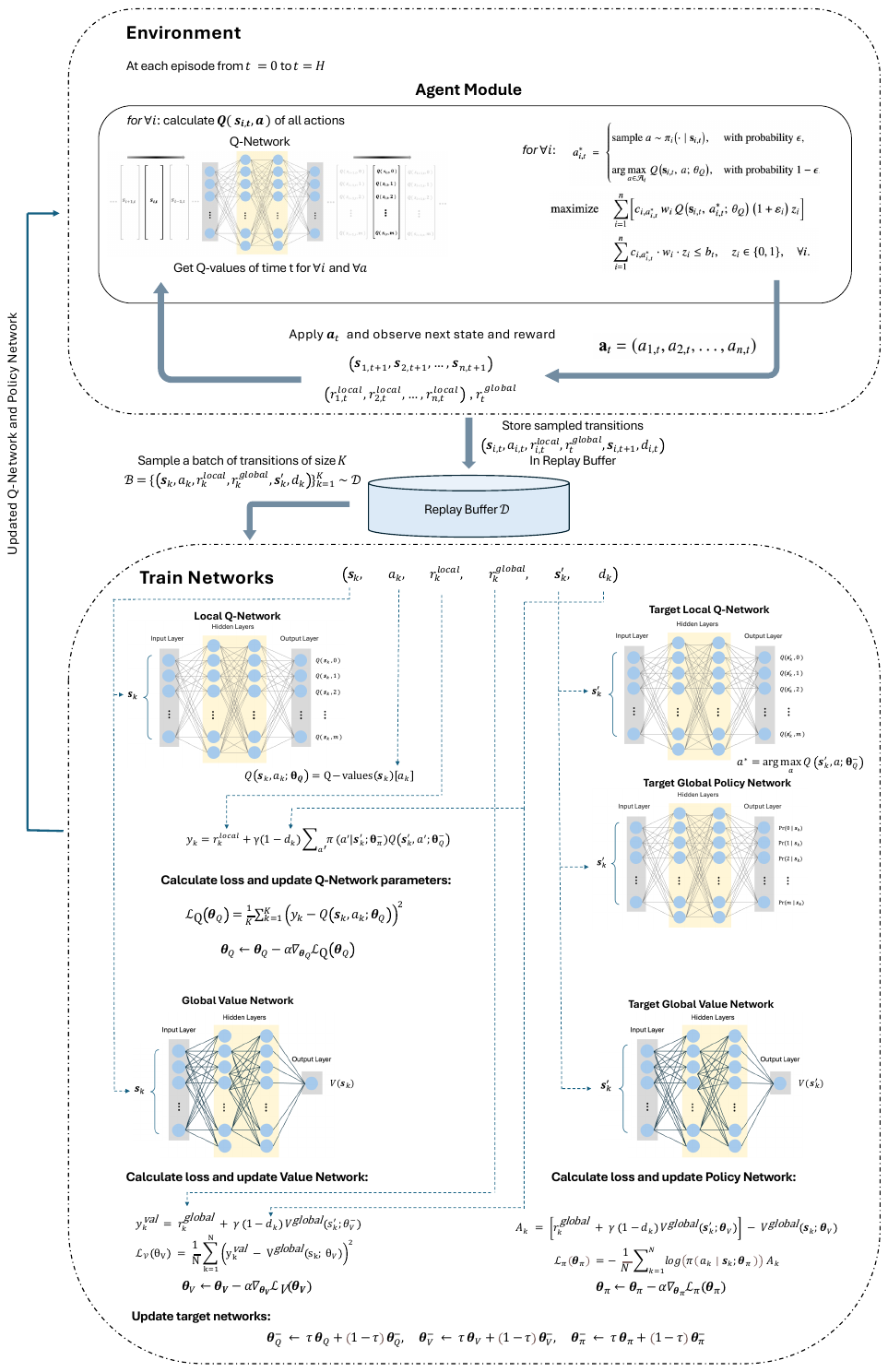}
    \caption{Flowchart of the proposed method}
    \label{fig:training_flowchart}
\end{figure}

The proposed method aims to overcome this over estimation bias with a novel decomposition framework. To achieve this, we first introduce a \emph{cost-normalized} reward at the local or asset level that expresses each intervention in terms of improvement per unit cost that makes Q-values more comparable across different assets in a heteregenous network.  Once every asset proposes its best local intervention, a network-wide resource allocation step enforces the annual budget constraint. This linear or knapsack-based optimization reconciles the local cost-effectiveness of each candidate action with the network-level budget constraint for that year. After that, we introduce an \emph{Expected SARSA} update that replaces the usual \(\max\) operator to estimate how likely each subsequent action will be chosen under a global policy which is governed by the annual budget constraint. By weighting next-step Q-values through the policy distribution, it can be determined which interventions can be afforded under a realistic multi-year budget. The separation of local Q-learning from budget enforcement maintains tractability, as each asset’s action-value function can be learned and updated independently, while the global resource allocation ensures feasible joint actions at each time step. A global policy and a global value network further guide exploration and learning toward actions that enhance overall network performance, ensuring that local decisions collectively yield long-term gains in LoS under budgetary constraints. The overall framework of the proposed method is illustrated in Fig.~\ref{fig:training_flowchart}. In the following, details of the key technical elements of the proposed method are first explained before the simulation and training of the networks are introduced.

\subsection{Cost-Normalized Reward and Local Q-Learning}

A key insight in this framework is to define each asset’s local reward in a \emph{cost-normalized} way, ensuring that intervention benefits are directly comparable across assets that may vary in size, importance, or unit costs. Specifically, if an action \(a_{i,t}\) is applied to asset \(i\) at time \(t\), let \(f_{i,a_{i,t},t}(s_{i,t})\) be the asset’s post-intervention condition, and let \(f_{i,0,t}(s_{i,t})\) be its condition if no intervention were undertaken. Denote the unit cost of action \(a_{i,t}\) by \(c_{i,a_{i,t},t}\), and the asset’s weight by \(w_i\). Then the immediate local reward is defined as
\begin{equation}
\label{eq:cost_normalized_reward}
  R_i\bigl(s_{i,t}, a_{i,t}\bigr)
  \;=\;
  \begin{cases}
    \displaystyle
    \frac{
      w_i\bigl[f_{i,a_{i,t},t}(s_{i,t}) - f_{i,0,t}(s_{i,t})\bigr]
    }{
      w_i\,c_{i,a_{i,t},t}
    }
    \;=\;
    \frac{
      f_{i,a_{i,t},t}(s_{i,t}) - f_{i,0,t}(s_{i,t})
    }{
      c_{i,a_{i,t},t}
    },
    & \text{if } a_{i,t} \neq 0,\\[8pt]
    0,
    & \text{if } a_{i,t} = 0.
  \end{cases}
\end{equation}
Since \(w_i\) cancels out in the numerator and denominator, \eqref{eq:cost_normalized_reward} is equivalently interpreted as the \emph{improvement in condition per unit cost}. Intuitively, the numerator quantifies how much a chosen intervention improves the network’s Level of Service (LoS) beyond natural deterioration, while the denominator is the intervention’s cost. By basing the local reward on this ratio, each asset’s Q-function captures how effective a given action is in “buying” condition improvements. 

To evaluate each asset’s actions over the long term, we introduce a Q-function, \(
Q_i\bigl(s_{i,t},\,a_{i,t}\bigr),
\)
that approximates the expected future return, expressed in these cost-normalized terms, from choosing action \(a_{i,t}\) in state \(s_{i,t}\) and subsequently following the learned policy. During each simulation step, a transition tuple 
\(\bigl(s_{i,t},\,a_{i,t},\,r_{i,t}^\textit{local},\,s_{i,t+1}\bigr)\) 
is observed, where \(r_{i,t}^\textit{local}\) is given by \eqref{eq:cost_normalized_reward}. Rather than using a max operator on the next-state actions, we adopt an \emph{Expected SARSA} approach, which integrates the policy’s probabilities over those actions. Specifically, let 
\(\pi_i\bigl(a' \mid s_{i,t+1}\bigr)\) 
be the probability of taking action \(a'\) in the next state \(s_{i,t+1}\). The Bellman equation for updating action value function then becomes:
\begin{equation}
\label{eq:expected_sarsa_target}
  Q_i\bigl(s_{i,t},\,a\bigr)
  \;=\;
  r_{i,t}^\textit{local}
  \;+\;
  \gamma
  \sum_{a' \,\in\, \mathcal{A}_i}
  \pi_i\bigl(a' \mid s_{i,t+1}\bigr)\,
  Q_i\bigl(s_{i,t+1},\,a'\bigr).
\end{equation}
This formulation lets the policy’s stochastic action selection influence how future returns are estimated. Actions deemed more likely under the current policy contribute more strongly to next-step Q-values, ensuring that exploration and exploitation are coordinated.

By iterating over such updates, each asset progressively refines \(Q_i(s_{i,t},\,a_{i,t})\) until it converges to an approximate \emph{cost-normalized} measure of the long-run benefit of each feasible intervention. Because rewards are scaled as “improvement per unit cost,” these local Q-values become directly comparable across assets, even when costs differ significantly.

Although each asset focuses on learning Q-values that reflect cost-effectiveness from its local perspective, a global resource allocation step ultimately decides which actions receive funding each year. This separation of concerns, involving local Q-learning for per-unit-cost effectiveness, combined with a budget allocation layer, maintains computational tractability and ensures that spending constraints are respected at the network level. As a result, each asset’s local improvement estimates are reconciled within a global optimization stage that maximizes overall Level of Service while adhering to annual budget limits.

\subsection{Budget Allocation via Linear Program}
\label{subsec:budget_lp}

While each asset independently learns a local Q-function that estimates the cost-normalized benefit of different actions, these candidate interventions must ultimately compete for funding under an annual budget~\(b_t\). To reconcile these local preferences with a global constraint, we solve a linear program (LP) each year to select which interventions receive funding. 

Define \(x_{i,j,t}\) as a binary decision variable that indicates whether action \(j\) is indeed executed for asset \(i\) (i.e., \(x_{i,j,t}=1\)) or not (\(x_{i,j,t}=0\)).  If a candidate action is chosen for asset~\(i\), its cost \(c_{i,j,t}\cdot w_i\) counts against the annual budget.  The objective is to pick a feasible subset of these interventions, maximizing total benefit.

Formally, suppose each local Q-function, \(Q_i^*\bigl(s_{i,t},\,j\bigr)\), represents the asset’s estimated \emph{improvement per unit cost} for action \(j\).  To translate these local Q-values into a consistent global objective, we multiply by \(c_{i,j,t}\,w_i\), giving an approximate total improvement at the network level.  The LP thus becomes:
\begin{align}
\label{eq:lp_objective}
\text{maximize} \quad
& 
\sum_{i=1}^{n} \sum_{j=0}^{m} 
  w_i\,c_{i,j,t}\,Q_i^*\bigl(s_{i,t},\,j\bigr) x_{i,j,t},
\\[5pt]
\label{eq:lp_budget}
\text{subject to} \quad
& 
\sum_{i=1}^{n} 
\sum_{j=1}^{m}
 c_{i,j,t} \, w_i \,x_{i,j,t} 
\;\;\le\; b_t,
\\[3pt]
\label{eq:lp_one_action}
& 
\sum_{j=0}^{m} x_{i,j,t} \;=\; 1,
\quad
x_{i,j,t}\in\{0,1\},
\quad
\forall i.
\end{align}
Constraint \eqref{eq:lp_budget} enforces that the total cost of selected actions does not exceed \(b_t\), while \eqref{eq:lp_one_action} ensures that exactly one action (including “No Intervention”) is chosen per asset.  The variables \(x_{i,j,t}=1\) if action \(j\) is selected for asset \(i\) at time~\(t\), and 0 otherwise.

Solving \eqref{eq:lp_objective}--\eqref{eq:lp_one_action} produces an integer solution \(\{x_{i,j,t}\}\) that identifies which interventions are funded under the budget.  Each asset’s final action \(a_{i,t}\) is then set to \(j\) if \(x_{i,j,t}=1\).  Because local Q-values already encode cost-normalized benefit, this LP offers a globally consistent way to allocate resources, ensuring that more cost-effective actions (as gauged by their learned Q-values) are favored, yet still respecting the shared budget constraint.

Overall, formulating budget allocation as an LP (or an efficient approximation) keeps the computational process tractable while maintaining a principled link between local Q-values and global resource constraints.   This synergy underlies the proposed method’s ability to scale to large networks without discarding budget feasibility or the advantages of reinforcement learning in discovering cost-effective maintenance actions. In addition, since this isolated, annual LP is a typical knapsacking problem, it can be quite effectively solved by a greedy algorithm using simple ranking.


\subsection{State Representation}

In the proposed DRL framework, rather than assigning a separate neural network to each asset, a unified (shared) neural network architecture is adopted. Asset-specific attributes are encoded as parametric inputs, allowing a single network model to handle many assets simultaneously. At each decision epoch \(t\), every asset \(i\) is described by a \emph{combined} or \emph{augmented} state vector, which includes both localized asset features and global, system-wide descriptors.

Locally, the asset state encompasses items such as the current condition indicator (e.g., a pavement quality index), relevant deterioration parameters, and cost information for potential maintenance actions. These details reflect the immediate, asset-specific reality of the pavement or infrastructure component, including its physical condition and resource requirements for each intervention option.

To inform each asset of the larger system context, the state vector is also augmented with global descriptors. By incorporating global information to the state vector, each asset learns to coordinate its decisions with overarching performance objectives and budget constraints, transcending purely local considerations.

Formally, the augmented state of asset \(i\) at time \(t\), denoted \(\mathbf{s}_{i,t}\), can be expressed as
\[
\mathbf{s}_{i,t}
\;=\;
\bigl(\,
\underbrace{s_{i,t},\,\lambda_i,\,c_{i,t},\,\dots}_{\text{local asset features}},\;\;
\underbrace{\overline{\text{LoS}}_t,\;\tfrac{b_r}{B},\;\tfrac{t}{h},\,\dots,\,\mathcal{H}_t}_{\text{global descriptors}}
\bigr),
\]
where \(s_{i,t}\) is the asset’s actual condition state at time \(t\), \(\lambda_i\) might be a deterioration rate parameter, \(c_{i,t}\) captures intervention cost information, \(\overline{\text{LoS}}_t\) is the network-wide performance at time \(t\), \(\tfrac{b_r}{B}\) indicates the fraction of budget remaining relative to a total budget \(B\), \(\tfrac{t}{h}\) is a temporal feature, and \(\mathcal{H}_t\) is a compressed distribution (e.g., a histogram) of all assets’ conditions.

The exact selection of features can be tailored to the specific application domain. In this study, the key principle is to incorporate vital local attributes alongside select global signals that allow each agent’s decisions to align with long-term, network-wide goals under finite budget constraints. In the notation used here, \(\mathbf{s}_{i,t}\) designates the augmented state for policy or Q-function updates, while \(s_{i,t}\) refers to the intrinsic condition state of the asset at time \(t\).

\subsection{Network Architecture}
\label{subsec:network_architecture}

\begin{figure}[H]
    \centering
    \begin{subfigure}[b]{0.45\textwidth}
        \centering
        \includegraphics[width=\textwidth]{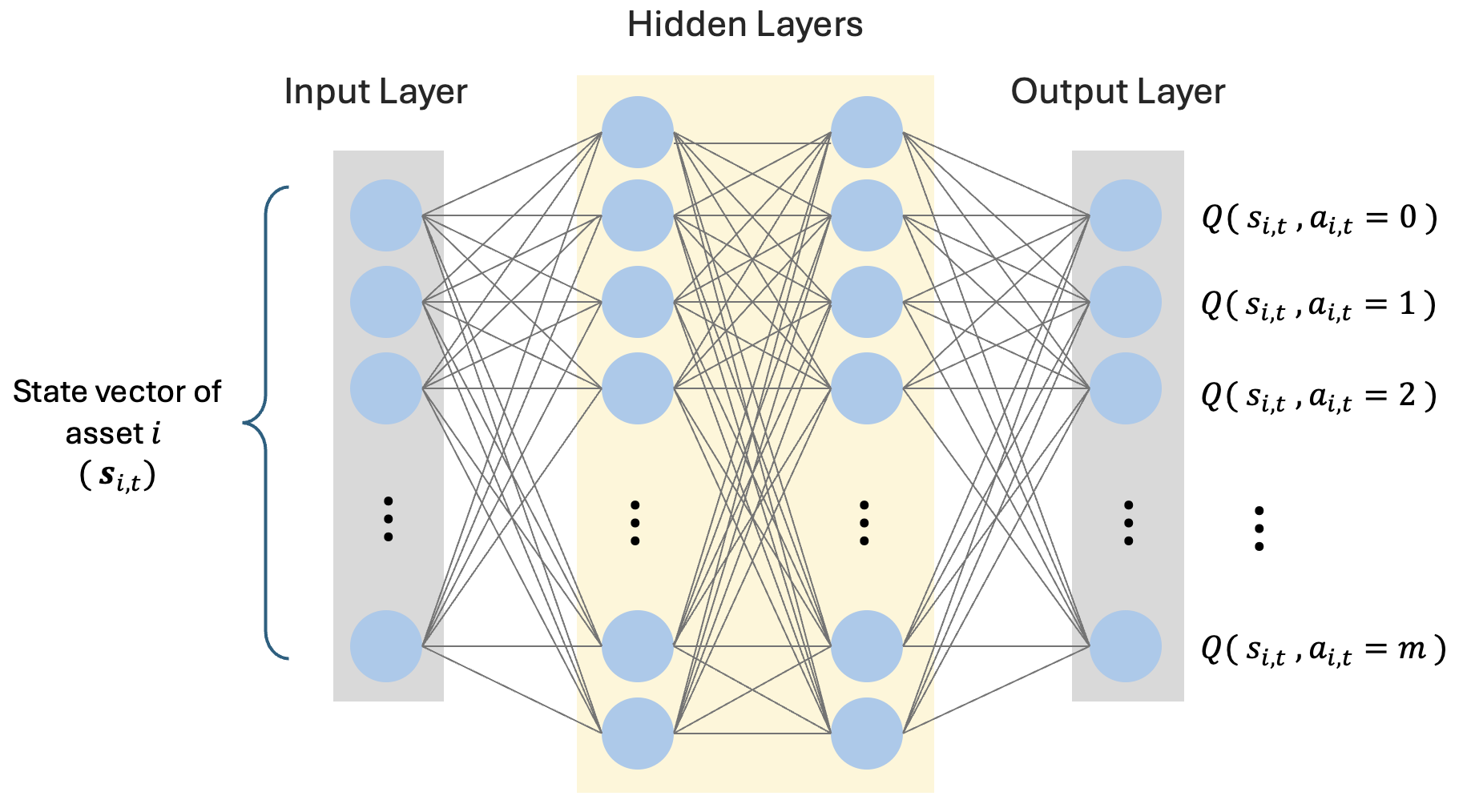}
        \caption{Q-Network Architecture}
        \label{fig:q_network}
    \end{subfigure}
    \hspace{0.05\textwidth}
    \begin{subfigure}[b]{0.45\textwidth}
        \centering
        \includegraphics[width=\textwidth]{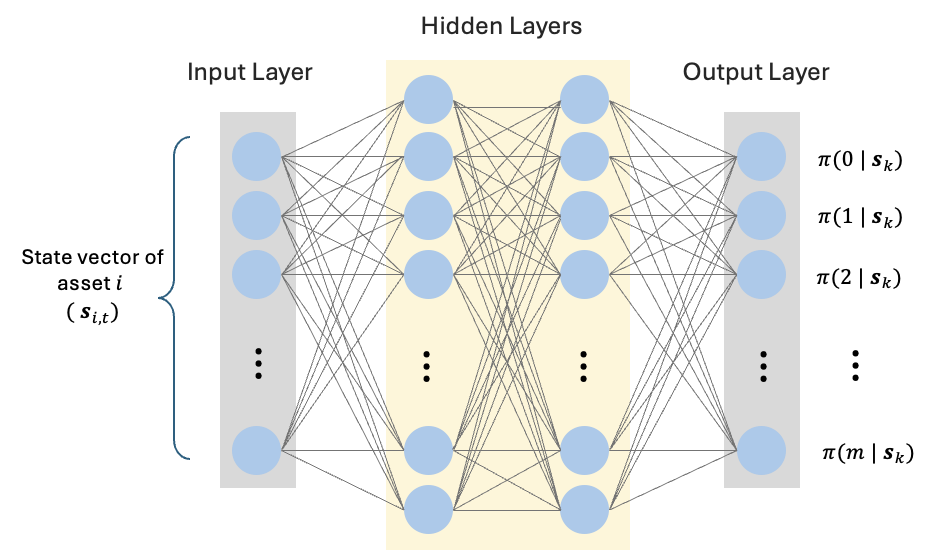}
        \caption{Policy Network Architecture}
        \label{fig:policy_network}
    \end{subfigure}

    \vspace{0.8em}  

    \begin{subfigure}[b]{0.45\textwidth}
        \centering
        \includegraphics[width=\textwidth]{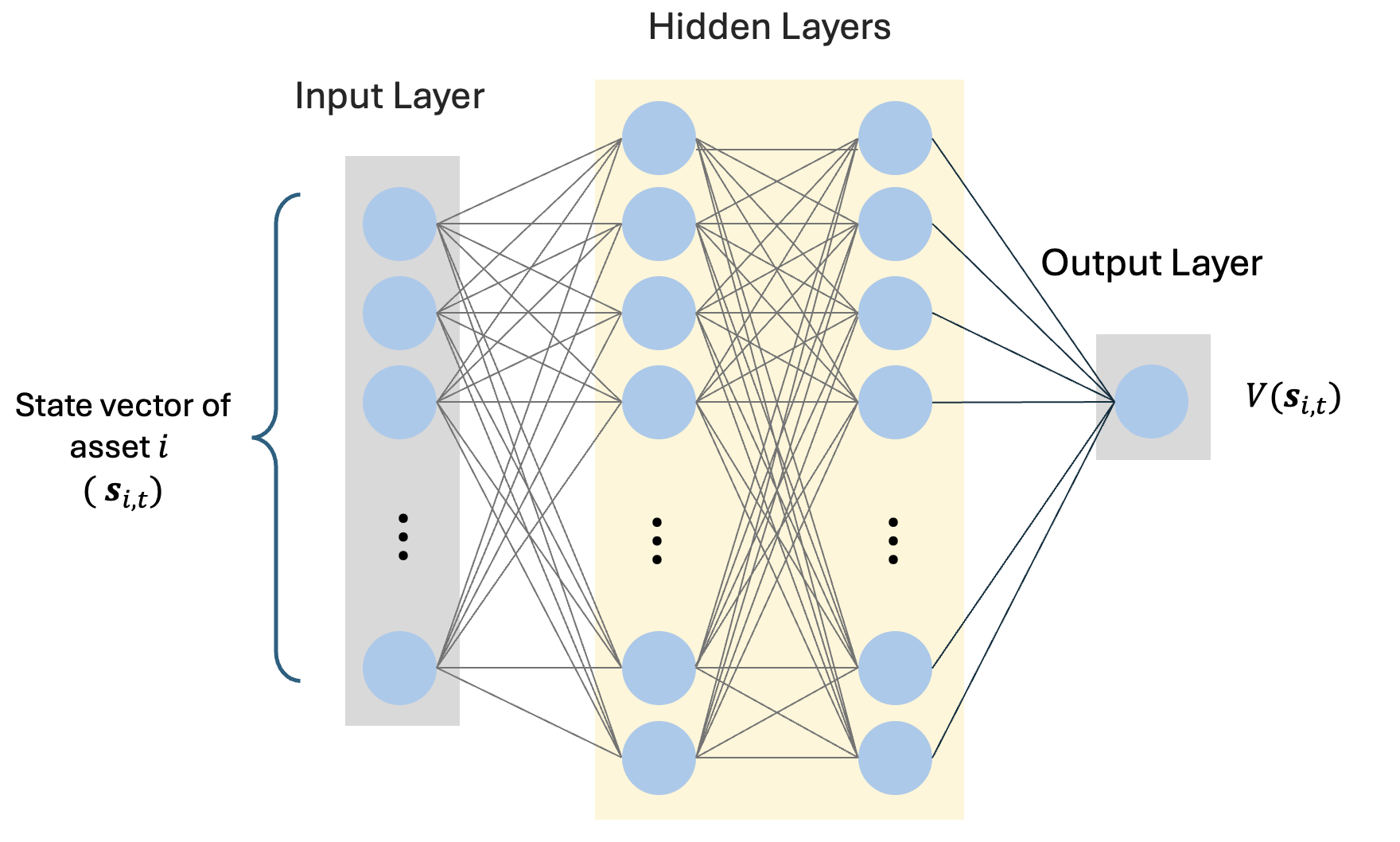}
        \caption{Value Network Architecture}
        \label{fig:value_network}
    \end{subfigure}

    \caption{Architectures of the local Q-network, global value network, and global policy network.}
    \label{fig:network_archs}
\end{figure}

To implement the proposed DRL framework, three neural networks are jointly employed: the \emph{Local Q Network}, the \emph{Policy Network}, and the \emph{Global Value Network}. Figure~\ref{fig:network_archs} illustrates the architecture of the neural networks. Each network shares a consistent input representation so that global context (e.g., budget or aggregate condition) can be combined with local asset features (e.g., current condition or deterioration rate). However, the three networks serve distinct roles in evaluating actions, shaping the policy, and estimating the system-wide return.

\vspace{0.5em}
\noindent
\textbf{Local Q Network:} 
This network estimates each asset’s \emph{cost-normalized} action-value function, \(Q_i(\mathbf{s}_{i,t},\,a_{i,t})\). Its inputs include the asset’s local state (condition, cost parameters, deterioration attributes) along with compressed network-level features (remaining budget, average LoS, time index). The output layer produces \(m+1\) values, one for each possible action, indicating the expected improvement per unit cost. By normalizing benefits by cost, the Q-values can be compared across different assets or interventions. During training, a TD update rule is applied, via Expected SARSA, which incorporates the probabilities from the Policy Network into the Q-value targets.

\vspace{0.5em}
\noindent
\textbf{Policy Network:}
Rather than simply selecting the action with the highest Q-value, the Policy Network learns a stochastic policy \(\pi_i(a_{i,t}\mid \mathbf{s}_{i,t})\) for each asset. Its inputs match those of the Local Q Network but it outputs a probability distribution over the \(m+1\) actions. This allows for exploration and for coordinating decisions with the global reward. Gradients from the \emph{global} objective (e.g., the network-wide LoS) flow back to shape these policy probabilities. Over time, the Policy Network balances taking apparently good interventions (exploitation) with sampling alternatives that might yield better outcomes (exploration).

\vspace{0.5em}
\noindent
\textbf{Global Value Network:}
This network provides a baseline for policy optimization by predicting the \emph{system-wide} future return from a given state. It ingests both asset-level state as well as network-level features summarizing the network (e.g., distribution of asset conditions, total budget usage) and outputs a single scalar, \(V^\text{global}(\mathbf{s}_{i,t})\). During policy updates, the \emph{advantage} is computed by comparing the observed global reward to the baseline predicted by this network. The result helps reduce variance when adjusting the Policy Network’s parameters: if the global reward surpasses the predicted baseline, the selected action is viewed favorably; otherwise, the policy is nudged away from that choice.

Although assets have distinct local states, the same neural architectures and parameter sets can be shared across assets, relying on input vectors to differentiate one asset from another. This approach improves scalability by avoiding an explosion in model size. When the parameter-sharing strategy is adopted, all assets use a single Local Q Network, a single Policy Network, and the same Global Value Network, thereby learning from the collective experience of the entire infrastructure system.

Overall, this three-network architecture links local decisions and global objectives in a unified DRL scheme. The Local Q Network furnishes cost-effectiveness estimates, the Policy Network ensures continuous exploration and adaptation toward higher global rewards, and the Global Value Network anchors policy updates against a system-level baseline. By sharing parameters across assets, the architecture scales to large networks, and by employing a budget allocation module, it preserves fidelity to real-world resource limitations.

\subsection{Simulation and Training the Networks}
\label{subsec:simulation_training}

This subsection describes how the unified Q-network, the global value network, and the policy network jointly learn an optimal multi-year maintenance strategy under budget constraints. The process unfolds iteratively across episodes, each simulating a full planning horizon of \(h\) years. Figure~\ref{fig:training_flowchart} provides a schematic overview.

\subsubsection{Action Selection and Budget-Constrained Optimization with Random Exploration}
At each decision epoch \(t \in \{1,\dots,h\}\), each asset \(i\) has an augmented local state \(\mathbf{s}_{i,t}\) and faces an action set \(\mathcal{A}_i\) (e.g., \(\{0,1,\dots,m\}\)). The local Q-network estimates the action-value function $ Q\bigl(\mathbf{s}_{i,t},\,a_i;\,\theta_Q\bigr)$, which measures the expected long-term improvement \emph{per unit cost} of selecting action \(a_i\). A naive approach would enumerate every \(\bigl(\mathbf{s}_{i,t},\,a\bigr)\) across all assets and solve a high-dimensional combinatorial program to satisfy the budget \(b_t\). However, such an approach becomes intractable for large-scale systems. Instead, a three-step action selection technique is proposed here. 
\begin{enumerate}
  \item [Step 1:] Determine Candidate Action per Asset. \\[3pt]
    Each asset \(i\) obtains a candidate action \(a_{i,t}^*\) through a hybrid of greedy choice and probabilistic exploration. With probability \(\epsilon\), the action is sampled from the learned policy network 
    \(\pi_i(\cdot \mid \mathbf{s}_{i,t})\), 
    while with probability \(1-\epsilon\), it is chosen by a greedy rule:
    \begin{equation}
      \label{eq:best_per_asset}
      a_{i,t}^*
      \;=\;
      \begin{cases}
        \displaystyle
        \text{sample } a \sim \pi_i\bigl(\cdot \mid \mathbf{s}_{i,t}\bigr),
        & \text{with probability }\epsilon,\\[6pt]
        \displaystyle
        \arg\max_{a \in \mathcal{A}_i}\,
        Q\bigl(\mathbf{s}_{i,t},\,a;\,\theta_Q\bigr),
        & \text{with probability }1-\epsilon.
      \end{cases}
    \end{equation}
    In early training, \(\epsilon\) is relatively high, encouraging extensive exploration via the stochastic policy \(\pi\). As learning converges, \(\epsilon\) typically decays, emphasizing exploitation of Q-value optima.

  \item [Step 2:] Network-Level Fund Allocation via LP. \\[3pt]
    Let \( w_i \, c_{i,a_{i,t}^*}\) be the cost of the candidate action \(a_{i,t}^*\) for asset \(i\). Deciding which candidate actions receive the current-year budget \(b_t\) is formulated as the following linear program (LP):
    \begin{equation}
      \label{eq:lp_obj_train}
      \text{maximize}
      \quad
      \sum_{i=1}^{n}
      \Bigl[
        w_i\, c_{i,a_{i,t}^*}\,
        Q\bigl(\mathbf{s}_{i,t},\,a_{i,t}^*;\,\theta_Q\bigr)\,
        \bigl(1 + \varepsilon_i\bigr)\,
        z_i
      \Bigr],
    \end{equation}
    subject to
    \begin{equation}
      \label{eq:lp_budget_train}
      \sum_{i=1}^{n}
      \Bigl[c_{i,a_{i,t}^*}\,w_i\Bigr]\,
      z_i
      \;\le\; b_t,
      \quad
      z_i \in \{0,1\},\quad
      i=1,\dots,n.
    \end{equation}
    Each binary variable \(z_i\) indicates whether asset \(i\) actually executes its candidate action (\(z_i=1\)) or defaults to no intervention (\(z_i=0\)). Multiplying the local Q-value by the action cost \(c_{i,a_{i,t}^*}\,w_i\) approximates the total expected improvement contributed by asset \(i\). The factor \(\bigl(1+\varepsilon_i\bigr)\), with \(\varepsilon_i \sim \mathcal{N}(0,\epsilon)\), injects additional exploration into the LP solver, permitting occasionally lower-Q actions to be selected.  A greedy knapsack-like heuristic is adopted: Specifically, assets’ candidate actions are ranked by 
\(
Q\bigl(\mathbf{s}_{i,t},\,a_{i,t}^*\bigr) (1+\varepsilon_i),
\)
and are selected in descending order until the budget is exhausted. Although not always perfectly optimal, this approach converges asymptotically to the LP solution for large-scale knapsack problems. Considering the random exploration included in the step, this computation-saving technique performs very well in practice.

  \item[Step 3:] Action Execution. \\[3pt]
    Based on the LP solution \(\{z_i\}\), if \(z_i=1\), then asset \(i\) applies the proposed action \(a_{i,t}^*\). Otherwise, it carries out the no-intervention option \(\,a_{i,t}=0\).
    Each asset’s state then evolves from \(s_{i,t}\) to \(s_{i,t+1}\) according to its deterioration and maintenance dynamics. 

    After deciding and executing \(\{a_{i,t}\}\) in year~\(t\), each asset \(i\) observes: $(\mathbf{s}_{i,t},\; a_{i,t},\; r_{i,t}^\textit{local}$, $ r_{i,t}^{\text{global}},\; \mathbf{s}_{i,t+1},\; d_{i,t} )$, where \(r_{i,t}^\textit{local}=R_i(s_{i,t},\,a_{i,t})\) is the local, cost-normalized reward and \(d_{i,t}\) is a terminal flag set if \(t=h\). These transitions are stored in a replay buffer that decouples the correlation of sequential observations. At the end of an episode (or periodically within it), mini-batches from this buffer are sampled to update the networks via stochastic gradient descent.
\end{enumerate}

Note that the framework employs two distinct exploration mechanisms: First, \(\epsilon\)-greedy choices in Step 1, where “random” means sampling from the stochastic policy \(\pi_i\). Second, a Gaussian factor \((1+\varepsilon_i)\) in the LP objective, letting candidate actions with slightly lower Q-values sometimes be funded. Early training uses larger \(\epsilon\) values, diversifying action selection. Over time, these exploration parameters decay, moving the agent toward deterministic strategies that leverage the best-known interventions to maximize network LoS under budget constraints.



\subsubsection{Loss Functions and Network Training}

To incorporate both local and global performance objectives, the training procedure updates three neural networks iteratively.
Let \(\theta_Q\), \(\theta_V\), and \(\theta_\pi\) denote the parameters of the local Q, global value, and policy networks, respectively. At each training step, a mini-batch of size \(N\) is drawn from the replay buffer: $ \bigl\{(\mathbf{s}_k,\;a_k,\;r_k^{\text{local}},\;r_k^{\text{global}},\;\mathbf{s}_k',\;d_k)\bigr\}_{k=1}^N $, where the subscript $k$ is used here to denote a sample randomly selected from the replay buffer. 
We assume a discount factor \(\gamma\in(0,1]\). Each network has a delayed (target) copy, \(\theta_Q^-\), \(\theta_V^-\), used to stabilize temporal-difference updates.

The local Q-value \(Q\bigl(s,a;\theta_Q\bigr)\) approximates the expected sum of discounted \emph{local} rewards \(\{r_{t}^{\text{local}}\}\). A standard temporal-difference (TD) target is constructed as
\[
y_k^{\text{TD}}
\;=\;
r_k^{\text{local}}
\;+\;
\gamma\,(1 - d_k)\,\max_{a'\,\in\,\mathcal{A}}Q\bigl(\mathbf{s}_k',\,a';\,\theta_Q^-\bigr),
\]
and the Q-network’s mean-squared error (MSE) loss is
\begin{equation}
\label{eq:LQ}
\mathcal{L}_Q(\theta_Q)
\;=\;
\frac{1}{N}\,
\sum_{k=1}^N
\Bigl(
y_k^{\text{TD}}
\;-\;
Q\bigl(\mathbf{s}_k,\,a_k;\,\theta_Q\bigr)
\Bigr)^2.
\end{equation}
Gradient descent on \(\nabla_{\theta_Q}\,\mathcal{L}_Q\) aligns the Q-network with observed local returns.

For the global value network,  a value function \(V^\text{global}\!(\mathbf{s};\theta_V)\) is trained using the global reward \(\{r_t^{\text{global}}\}\) to capture long-term \emph{network-wide} performance. For transition \(\bigl(\mathbf{s}_k,\,r_k^{\text{global}},\,\mathbf{s}_k'\bigr)\), a one-step bootstrap target is
\[
y_k^{\text{val}}
\;=\;
r_k^{\text{global}}
\;+\;
\gamma\,(1 - d_k)\,
V^\text{global}\!\bigl(\mathbf{s}_k';\theta_V^-\bigr),
\]
and the MSE loss is
\begin{equation}
\label{eq:LV}
\mathcal{L}_V(\theta_V)
\;=\;
\frac{1}{N}\,\sum_{k=1}^N
\Bigl(
y_k^{\text{val}}
\;-\;
V^\text{global}\!\bigl(\mathbf{s}_k;\,\theta_V\bigr)
\Bigr)^2.
\end{equation}
Minimizing \(\mathcal{L}_V\) trains the global value network to predict future LoS-based returns from each state.

The policy network \(\pi(a\mid \mathbf{s};\theta_\pi)\) outputs a probability distribution over actions, combining local and global signals. Its updates use a policy gradient with an advantage function derived from \(V^\text{global}\). For transition \(k\), the advantage is
\[
A_k
\;=\;
\bigl[r_k^{\text{global}}
\;+\;
\gamma\,(1 - d_k)\,
V^\text{global}\!\bigl(\mathbf{s}_k';\theta_V\bigr)\bigr]
\;-\;
V^\text{global}\!\bigl(\mathbf{s}_k;\theta_V\bigr).
\]
The policy loss is
\begin{equation}
\label{eq:Lpi}
\mathcal{L}_\pi(\theta_\pi)
\;=\;
-\;\frac{1}{N}
\,\sum_{k=1}^{N}
\log\Bigl(\pi\bigl(a_k\mid \mathbf{s}_k;\theta_\pi\bigr)\Bigr)\,
A_k.
\end{equation}
By descending \(\nabla_{\theta_\pi}\,\mathcal{L}_\pi\), actions that yield higher-than-expected global returns are increasingly favored.

At each iteration, we perform parameter updating using the following:
\[
\theta_Q
\;\leftarrow\;
\theta_Q - \alpha_Q\,\nabla_{\theta_Q}\,\mathcal{L}_Q(\theta_Q),\quad
\theta_V
\;\leftarrow\; 
\theta_V - \alpha_V\,\nabla_{\theta_V}\,\mathcal{L}_V(\theta_V),\quad
\theta_\pi
\;\leftarrow\;
\theta_\pi - \alpha_\pi\,\nabla_{\theta_\pi}\,\mathcal{L}_\pi(\theta_\pi),
\]
where \(\{\alpha_Q,\alpha_V,\alpha_\pi\}\) are learning rates. Periodically, each target network is softly updated, e.g.
\[
\theta_Q^-\;\leftarrow\;\tau\,\theta_Q + (1-\tau)\,\theta_Q^-,\quad
\theta_V^-\;\leftarrow\;\tau\,\theta_V + (1-\tau)\,\theta_V^-,\quad
\theta_\pi^-\;\leftarrow\;\tau\,\theta_\pi + (1-\tau)\,\theta_\pi^-,
\]
to stabilize temporal-difference targets.

\section{Results and Discussion}
\subsection{Case Study}

This case study focuses on a large-scale pavement network containing 68{,}800 individual pavement segments distributed throughout a metropolitan area. The entire network spans a total pavement area of approximately 59{,}856{,}743.2 square meters, which is equivalent to 17{,}101.9 lane-kilometres at a 3.5\,m lane width. Each pavement segment features distinct characteristics, such as size, current condition, deterioration profile, and maintenance costs, reflecting the variety typically observed in real-world infrastructure networks.

The objective is to determine an optimal 20-year maintenance plan that maximizes the network’s Level of Service (LoS) while respecting a \$200\,million annual budget. Pavement condition is measured by the Pavement Quality Index (PQI), a composite indicator that combines factors like surface smoothness, rutting, and surface distress. The overall network LoS is then represented by a weighted average of segment-level PQI values. Specifically, if $\text{PQI}_{i,t}$ denotes the PQI for segment \(i\) in year \(t\), and $w_i$ is its associated weight (e.g., area or traffic importance), then the network LoS at year \(t\) is
\begin{equation}
  \overline{\text{PQI}}_t 
  \;=\;
  \frac{\sum_{i=1}^{n} \,w_i \,\text{PQI}_{i,t}}{\sum_{i=1}^{n} \,w_i}.
\end{equation}

The network comprises three classes of roads: Arterials, Collectors, and Locals. Arterials constitute 30.4\% of the total pavement area, Collectors 19.9\%, and Locals 49.7\%. Maintenance costs scale with road classification. Rehabilitation, for instance, ranges from \$40 (Arterial) to \$20 (Local) per square meter, whereas reconstruction ranges from \$200 (Arterial) down to \$150 (Local) per square meter.

\begin{figure}[ht]
  \centering
  \begin{subfigure}{0.48\textwidth}
    \centering
    \includegraphics[width=\textwidth]{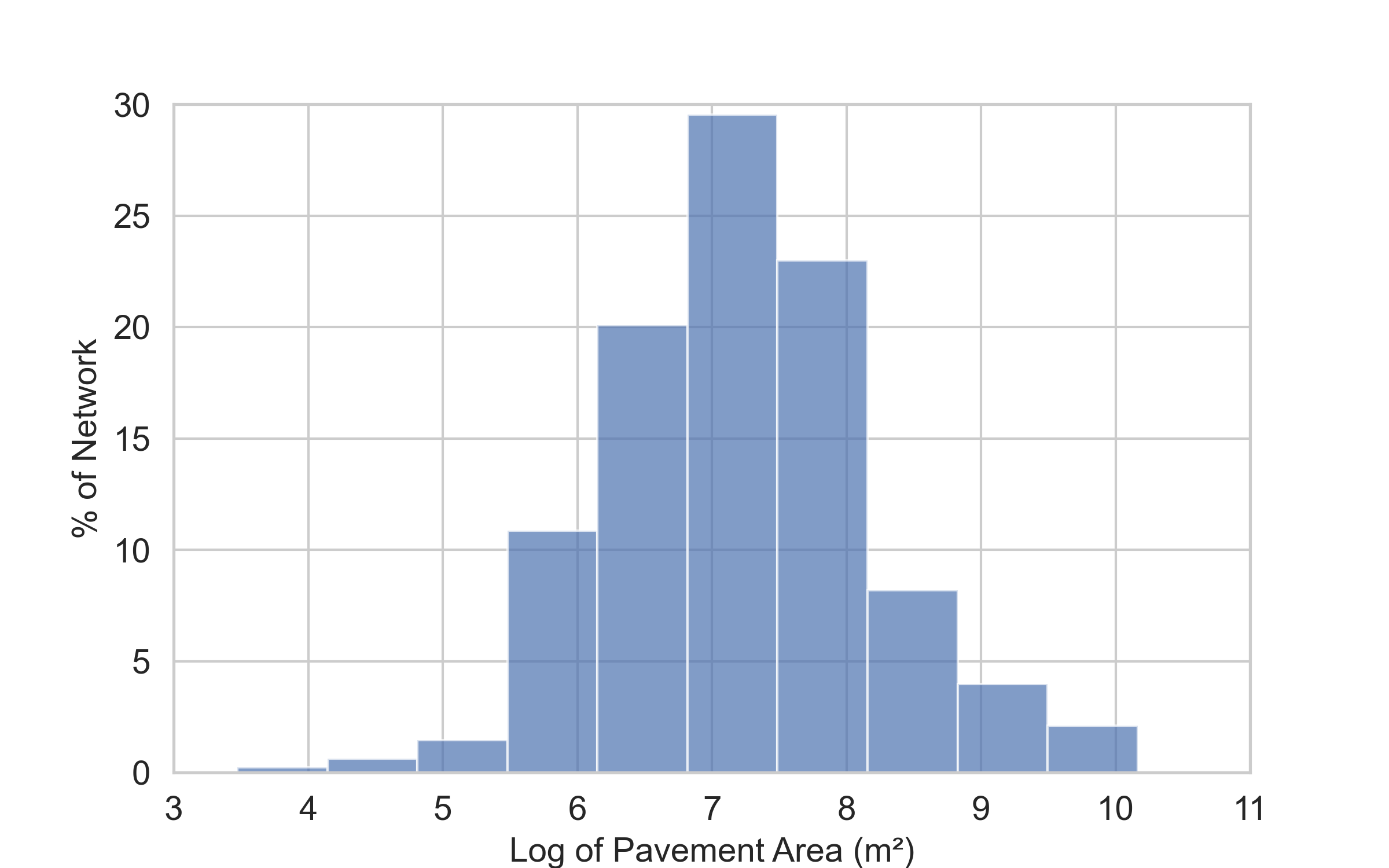}
    \caption{Pavement area}
    \label{fig:distribution_area}
  \end{subfigure}
  \hfill
  \begin{subfigure}{0.48\textwidth}
    \centering
    \includegraphics[width=\textwidth]{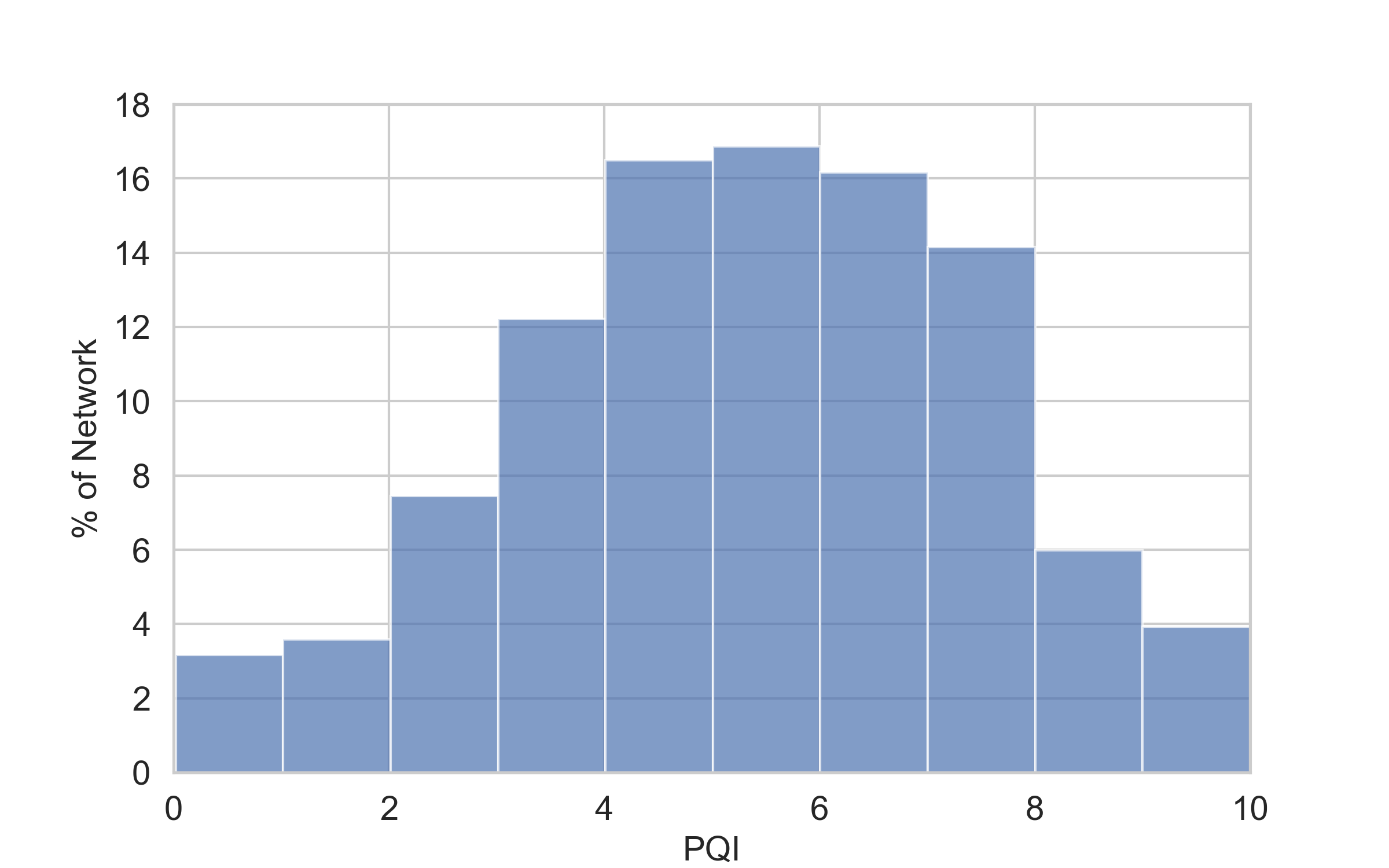}
    \caption{Initial PQI}
    \label{fig:initial_pqi_histogram}
  \end{subfigure}
  \caption{Histograms of pavement area and initial pavement condition across the network.}
  \label{fig:area_condition}
\end{figure}

Figure~\ref{fig:area_condition} illustrates the distribution of pavement area and the initial PQI for all segments. In Figure~\ref{fig:distribution_area}, the histogram reveals the variation in segment sizes, while Figure~\ref{fig:initial_pqi_histogram} presents the distribution of initial conditions observed across the network.

\vspace{1em}

\subsubsection{Performance Deterioration Model}
\begin{figure}[ht]
  \centering
  \begin{subfigure}{0.48\textwidth}
    \centering
    \includegraphics[width=\textwidth]{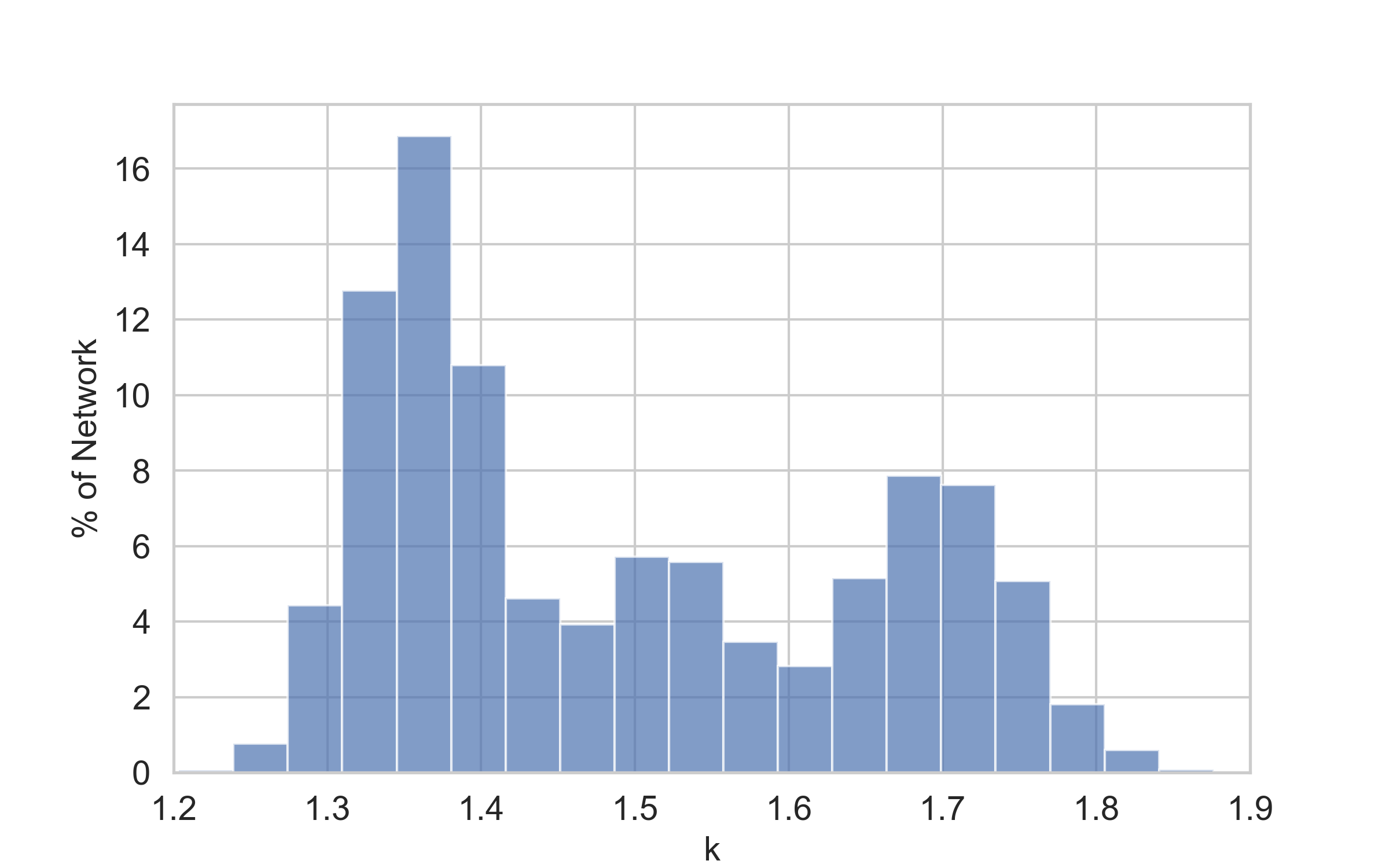}
    \caption{Shape parameter \(k\).}
    \label{fig:distribution_k}
  \end{subfigure}
  \hfill
  \begin{subfigure}{0.48\textwidth}
    \centering
    \includegraphics[width=\textwidth]{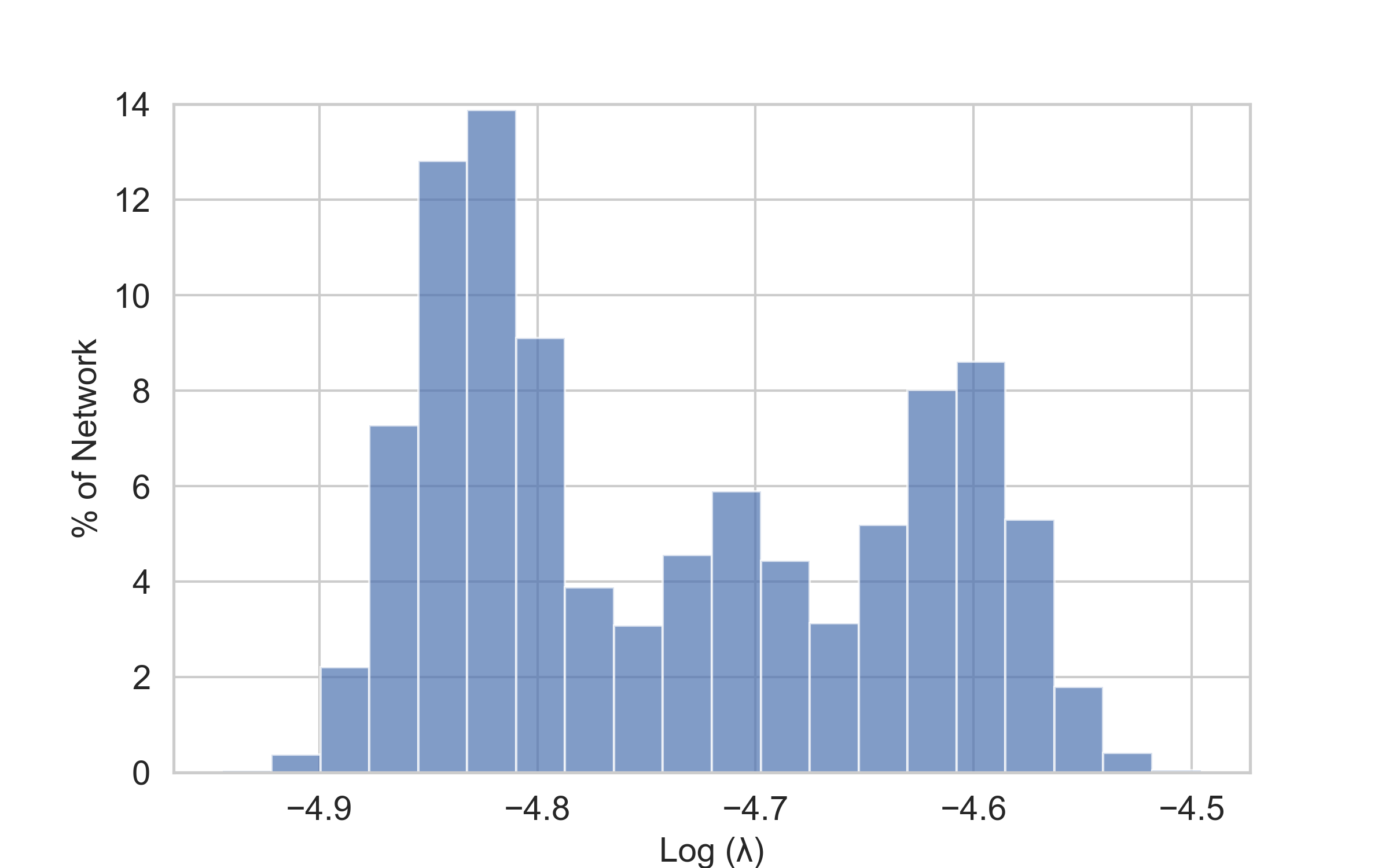}
    \caption{Scale parameter \(\lambda\).}
    \label{fig:distribution_lambda}
  \end{subfigure}
  \caption{Histograms of the Weibull parameters \(k\) and \(\lambda\) across the pavement network.}
  \label{fig:distribution_parameters}
\end{figure}

Each pavement segment deteriorates over time following a Weibull-type function:
\begin{equation}
  \text{PQI}_{i,t_{\text{age}}}
  \;=\;
  \text{PQI}_{\text{max}}\,
  \exp\bigl(-\,\lambda_i\,\tau^{k_i}\bigr),
\end{equation}
where $\text{PQI}_{\text{max}} = 10$ is the highest achievable PQI, $\lambda_i$ and $k_i$ are the deterioration parameters for segment \(i\), and $\tau$ indicates the time in years since the last major intervention. An incremental version of this expression can be used to model year-to-year PQI reductions. Figure~\ref{fig:distribution_parameters} shows how the Weibull parameters $\lambda_i$ and $k_i$ vary among segments, accounting for diverse environmental and traffic conditions.

\subsubsection{Intervention Actions and Maintenance Effectiveness}

Three major interventions are considered in the case study. They are: Do Nothing (0), Rehabilitation (1), and Reconstruction (2). While Do Nothing incurs no cost and the segment continues to deteriorate naturally based on the above model, Reconstruction restores the segment to \(\text{PQI}_{\text{max}}=10\), thereby providing a long-lasting effect but at the highest cost among the available options. In contrast, Rehabilitation provides a moderate improvement in the PQI, especially for segments in fair condition. If the PQI is above a threshold, a fixed increase is applied up to a maximum of 9.5. For lower PQI values, the improvement is scaled. Formally, after Rehabilitation, the PQI is changed to 
\begin{equation}
  \text{PQI}_{i,t+1} 
  \;=\; 
  \min \, \!\Bigl(\text{PQI}_{i,t} \;+\; \Delta_{\text{rehab}}\,
  	\frac{\text{PQI}_{i,t}}{9.5}, \;9.5\Bigr),
\end{equation}
where $\Delta_{\text{rehab}}=2.5.$

\subsubsection{State Representation and Input Dimensions}

A segment’s state vector integrates local and network-level data, ensuring that decisions remain aligned with both segment-specific needs and system-wide priorities. This vector for segment \(i\) in year \(t\) comprises 18 components:
\begin{equation}
  s_{i,t} 
  \;=\; 
  \Bigl(\,
    \text{PQI}_{i,t}, \, w_i, \, \lambda_i, \, k_i, \,
    c_{i,1}, \, c_{i,2}, \,
    \tfrac{t}{T}, \,
    \tfrac{b_r}{B}, \,
    \overline{\text{PQI}}_t, \dots,
    \text{10 histogram bins}
  \Bigr).
\end{equation}

Here, $\text{PQI}_{i,t}$ is the current condition, $w_i$ is the pavement area weight, $(\lambda_i, k_i)$ represent Weibull deterioration parameters, $(c_{i,1}, c_{i,2})$ denote unit costs for rehabilitation and reconstruction, $\tfrac{t}{T}$ is the normalized time index, $\tfrac{b_r}{B}$ is the ratio of remaining to total budget, and $\overline{\text{PQI}}_t$ is the average network PQI in year \(t\). A binned histogram of PQI values portrays the overall condition distribution, allowing the agent to account for how other assets’ conditions might affect budget usage and action funding.

This design captures both local and global considerations, supporting policy decisions that best serve individual segment health while respecting the \$200\,million annual spending limit and pursuing an optimal level of service over the 20-year horizon.
\subsection{Training Results and Analysis}

\begin{figure}[h!]
    \centering
    \begin{subfigure}{0.48\textwidth}
        \includegraphics[width=\textwidth]{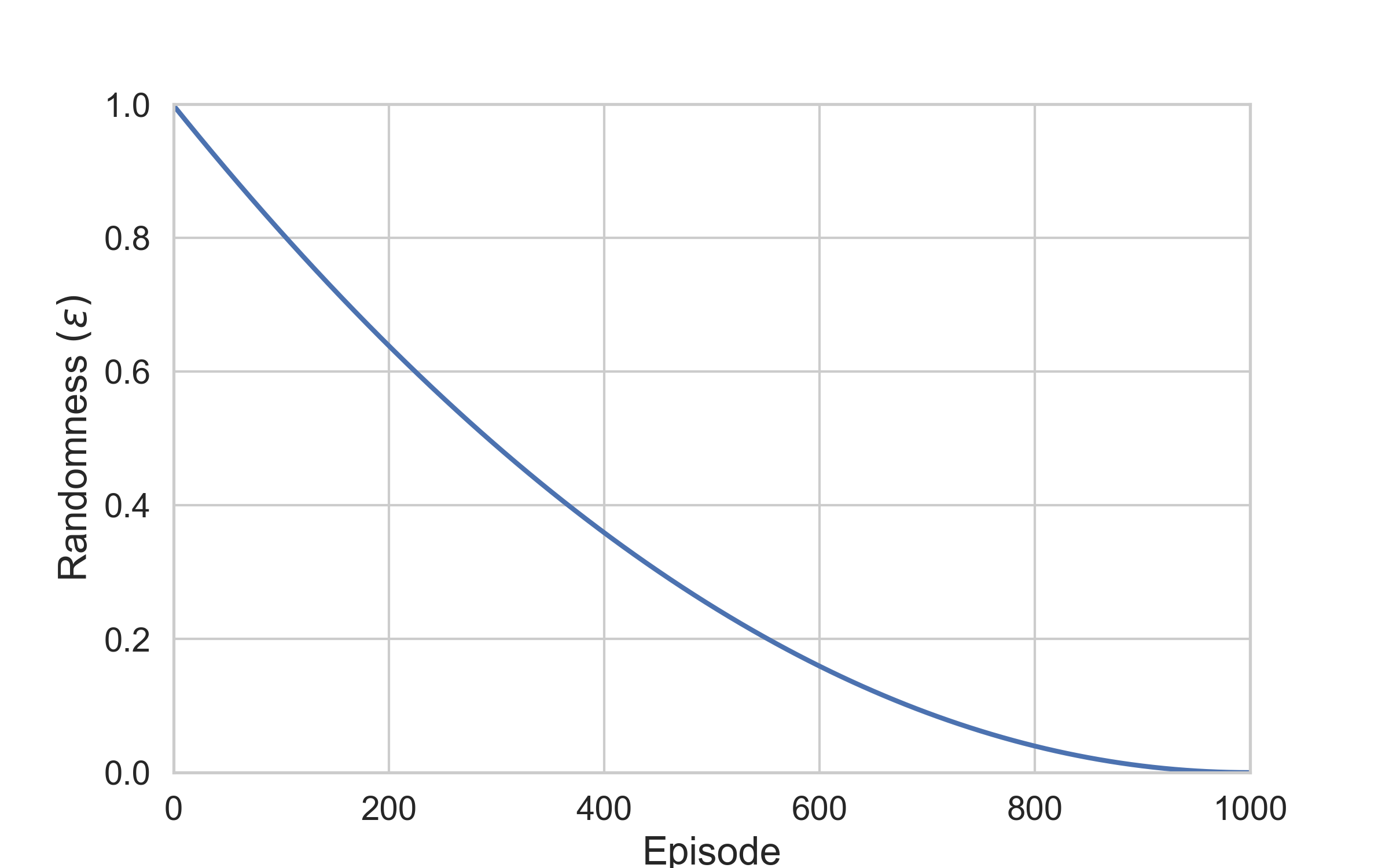}
        \caption{Randomness (\(\epsilon\))}
        \label{fig:randomness}
    \end{subfigure}
    \hfill
        \begin{subfigure}{0.48\textwidth}
        \includegraphics[width=\textwidth]{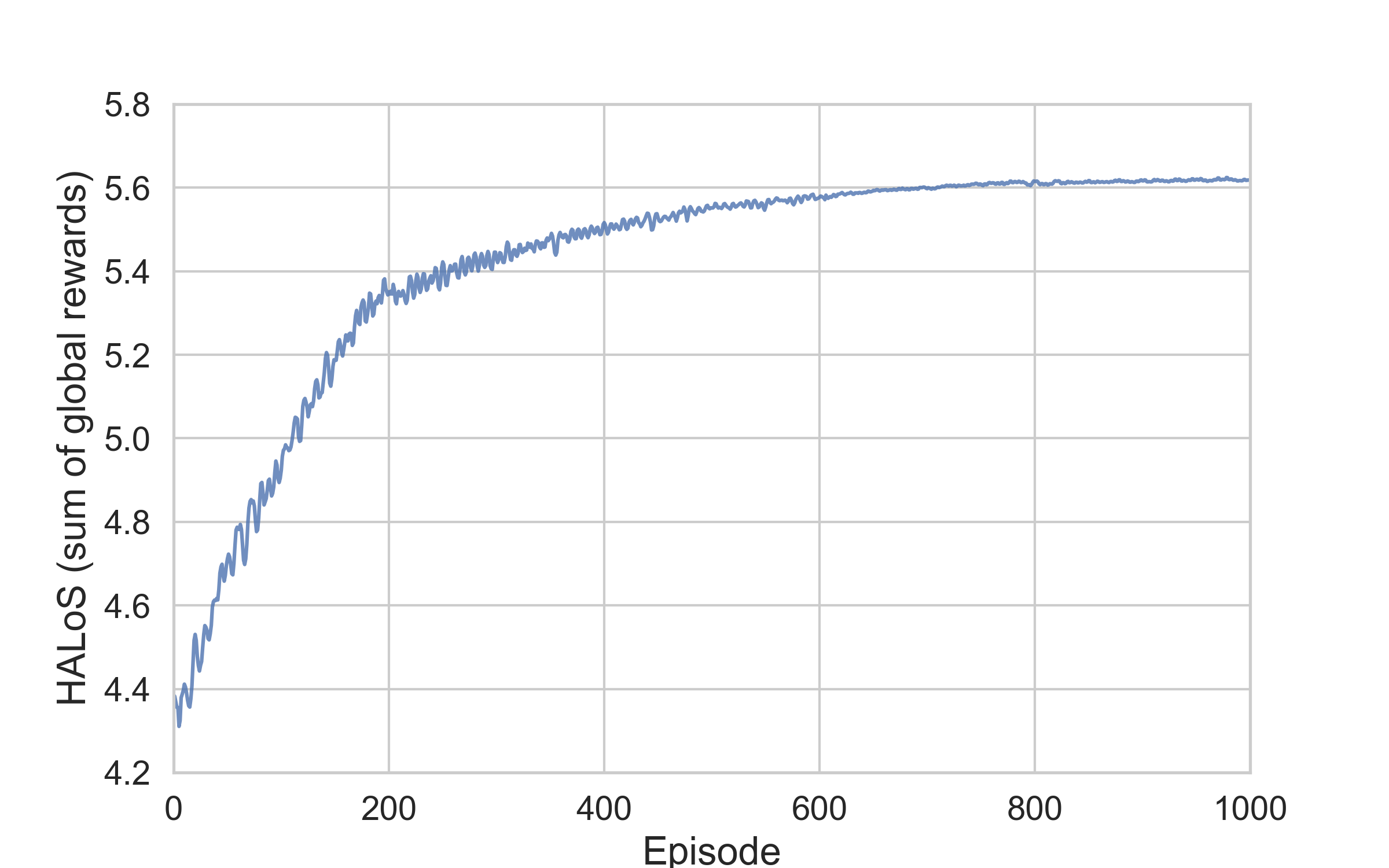}
        \caption{Average network LoS}
        \label{fig:reward}
    \end{subfigure}
    \vfill  
    \begin{subfigure}{0.48\textwidth}
        \includegraphics[width=\textwidth]{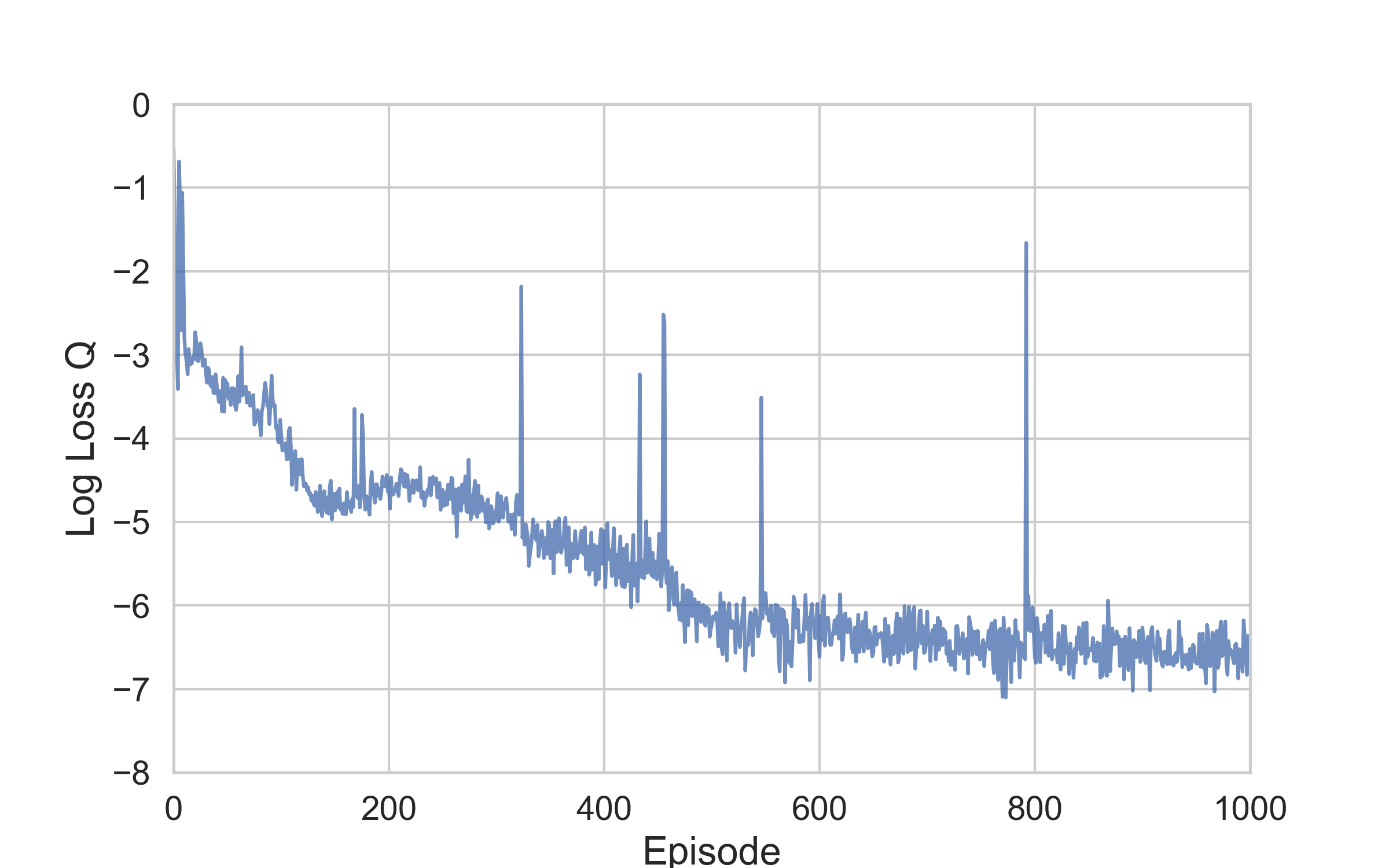}
        \caption{Local Q-network loss} 
        \label{fig:loss_q}
    \end{subfigure}
    \hfill
    \begin{subfigure}{0.48\textwidth}
        \includegraphics[width=\textwidth]{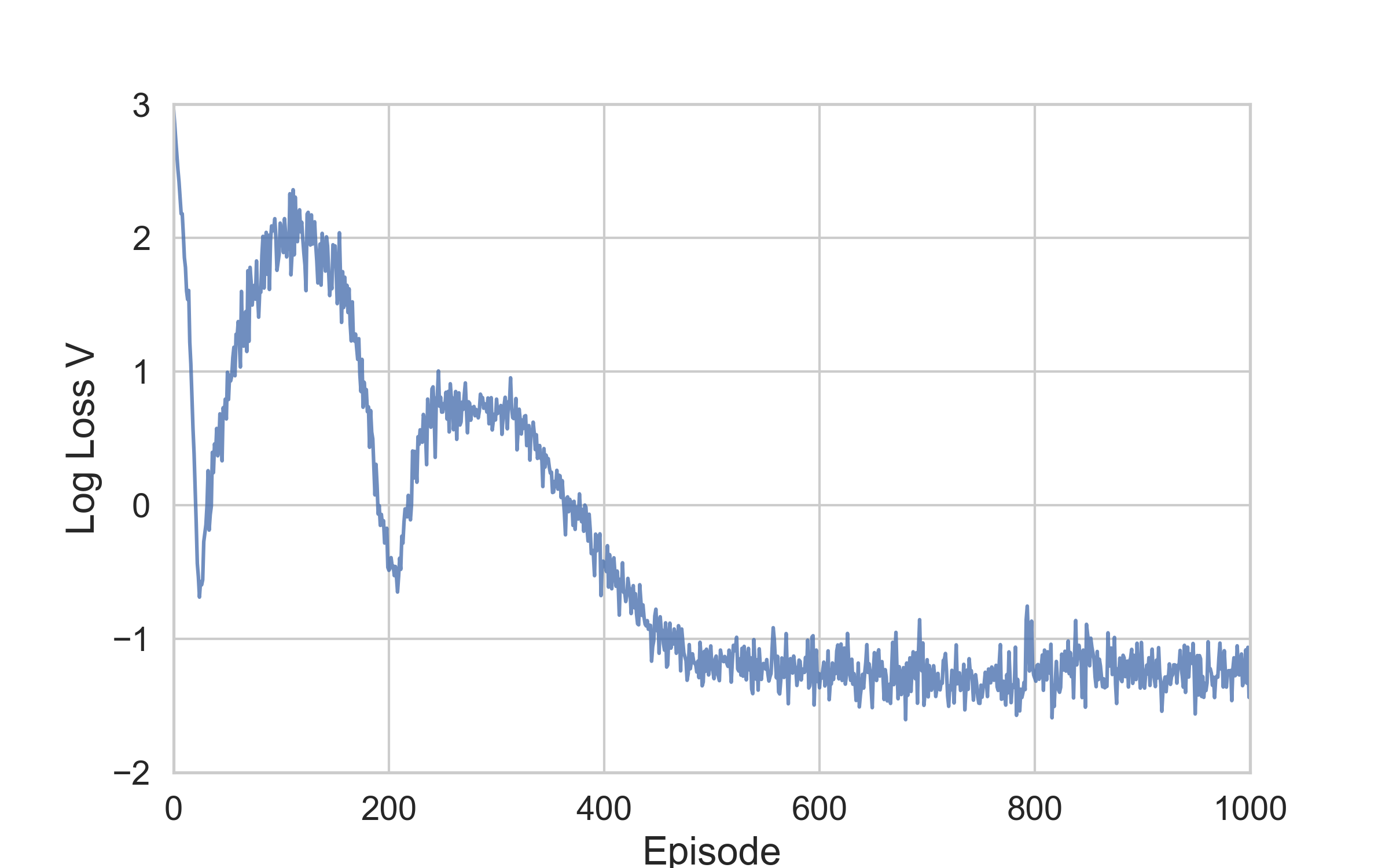}
        \caption{Global value network loss}
        \label{fig:loss_v}
    \end{subfigure}
    \caption{Trends of key training parameters over 1{,}000 episodes.}
    \label{fig:training}
\end{figure}

All experiments were implemented in \texttt{Python} using the \texttt{PyTorch} and \texttt{OR-tools} libraries for deep learning and LP solver, respectively. Training was performed on a workstation equipped with an Intel Core i9-13900K (3.00\,GHz) processor, 64\,GB of RAM, and no dedicated GPU usage. Under these conditions, running 1{,}000 training episodes required approximately 97\,minutes. A series of 1{,}000 training episodes was conducted to develop an optimal maintenance strategy aimed at maximizing the network’s horizon-averaged LoS, subject to annual budget constraints. During training, several performance indicators were tracked, including the average network LoS, the Q-network’s loss, the value network's loss, and the exploration randomness factor \(\epsilon\). These metrics capture the model’s learning dynamics and the evolution of its policy. Figure~\ref{fig:training} plots the trends of these quantities across the 1{,}000 episodes.

Exploration in the training process is governed by the factor \(\epsilon\), whose decay over episodes is shown in Figure~\ref{fig:randomness}. Early on, a higher degree of randomness fosters broad exploration of possible policies, preventing the agent from prematurely settling on suboptimal solutions. As training advances, \(\epsilon\) decreases, allowing the policy to focus on exploiting learned knowledge for improved, more stable decisions.

Figure~\ref{fig:reward} displays how the model’s cumulative reward evolves over time, where the reward corresponds to the sum of network LoS over each planning horizon in a single episode. This trajectory rises rapidly in the initial training phase, indicating that the model quickly discovers high-value maintenance actions. It then stabilizes, reflecting a learned policy that consistently maximizes the horizon-averaged LoS under budget constraints.

Figures~\ref{fig:loss_q} and \ref{fig:loss_v} illustrate the learning dynamics of both the local Q-network and the global value network. Initially, both losses tend to be relatively high as the networks adjust to approximate the value functions. Over time, these losses diminish, signaling that the framework is converging to a policy that effectively balances local actions with global LoS objectives.

\subsection{Performance Analysis}

A comparative evaluation was conducted to evaluate the proposed Network DQL method against three benchmarks: (1) a Worst-First strategy, (2) a progressive linear programming (LP) approach, and (3) a Hybrid LP-GA method developed by \citeN{Fard2024}. Table~\ref{tab:los_summary} presents the Horizon-Averaged LoS (HALoS) and the End-of-Horizon LoS (EHLoS) achieved by each method over a 20-year planning horizon.

\begin{table}[ht]
    \centering
    \caption{Comparison of the HALoS and EHLoS for the four solution methods.}
    \begin{tabular}{lcc}
        \hline
        \textbf{Method} & \textbf{HALoS} & \textbf{EHLoS} \\
        \hline
        Worst-First               & 4.55 & 4.00 \\
        Progressive LP            & 5.46 & 5.58 \\
        Hybrid LP-GA              & 5.50 & 5.68 \\
        Network DQL (proposed)    & 5.62 & 5.86 \\
        \hline
    \end{tabular}
    \label{tab:los_summary}
\end{table}

\begin{figure}[ht]
    \centering
    \includegraphics[width=0.85\textwidth]{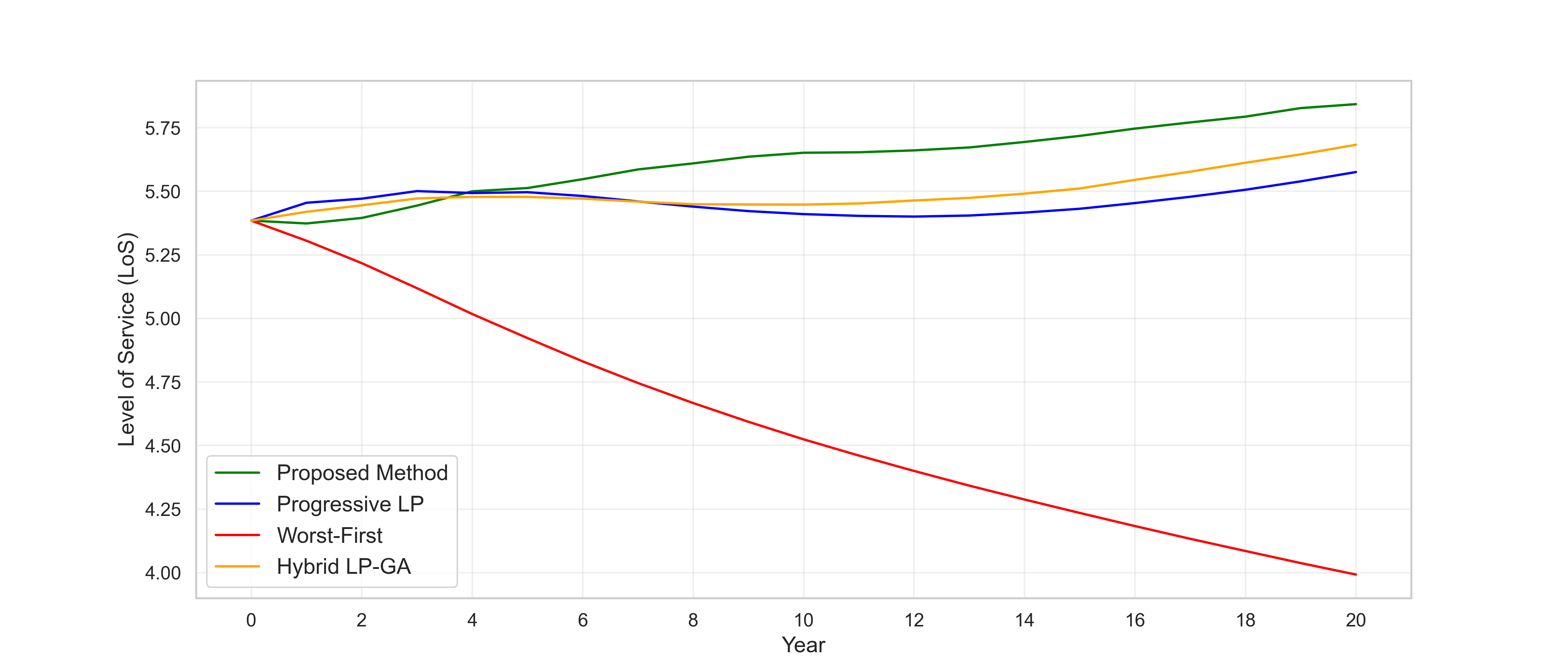}
    \caption{Average network LoS profiles for the four strategies (1,000 simulations).}
    \label{fig:comparison}
\end{figure}

Figure~\ref{fig:comparison} illustrates the LoS trajectories for these approaches, averaged over 1,000 simulations. The Worst-First strategy consistently allocates resources to the most deteriorated segments first, causing steep budget use on reconstruction every year. Although it temporarily addresses severely damaged segments, the LoS steadily decreases toward the horizon’s end, yielding a relatively low HALoS of 4.55. 

The Progressive LP method updates its plan each year by solving a binary linear program that maximizes next-year LoS under a fixed annual budget. This leads to an initial improvement but lacks long-term foresight, as evidenced by its lower LoS in later years. The Hybrid LP-GA approach performs better than Progressive LP by blending local search (LP) with more global search capabilities (GA). This mixed strategy attains a higher LoS than Progressive LP (see Table~\ref{tab:los_summary}), though it still falls short of the proposed method’s performance. In contrast, the Network DQL framework optimizes across the entire 20-year span and takes into account future payoffs from each maintenance action. This strategy achieves a HALoS of 5.62 and an EHLoS of 5.86, surpassing the other three methods. 

\begin{figure}[ht]
    \centering
    \begin{subfigure}{0.48\textwidth}
        \includegraphics[width=\textwidth]{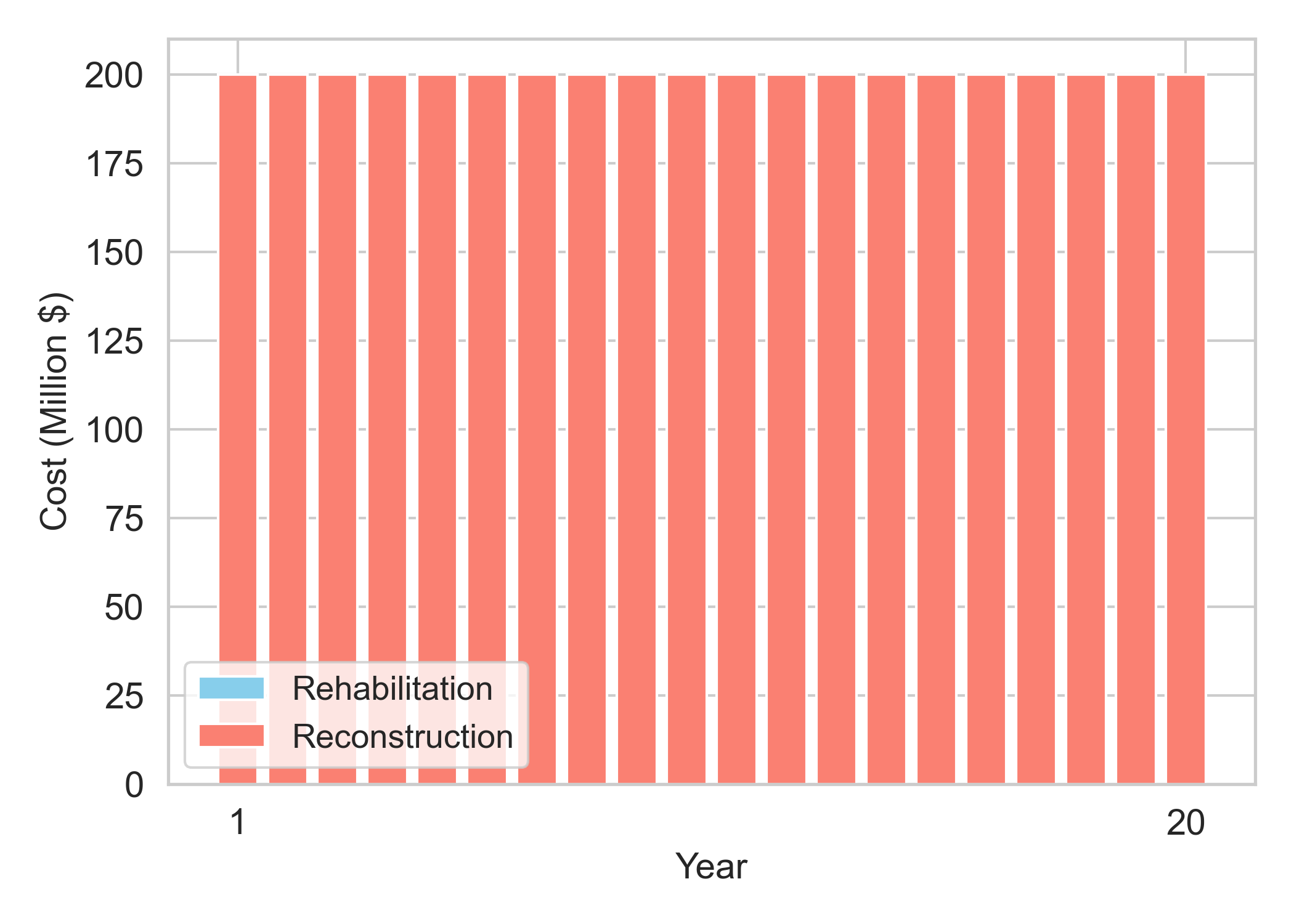}
        \caption{Worst-First strategy}
        \label{fig:worst_first_annual_cost}
    \end{subfigure}
    \hfill
    \begin{subfigure}{0.48\textwidth}
        \includegraphics[width=\textwidth]{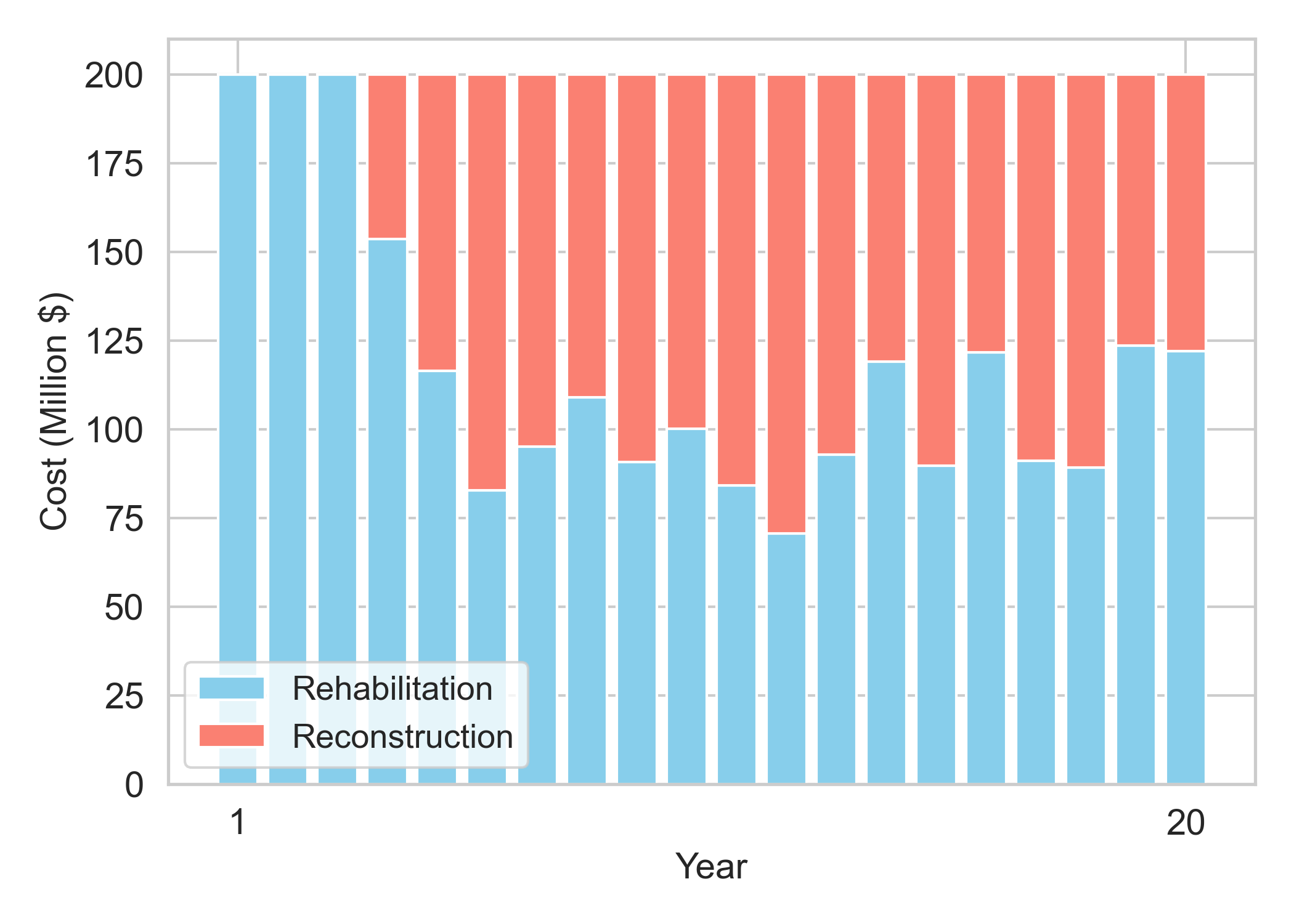}
        \caption{Progressive LP}
        \label{fig:prog_annual_cost}
    \end{subfigure}
    
    \vspace{0.5cm}  
    
    \begin{subfigure}{0.48\textwidth}
        \includegraphics[width=\textwidth]{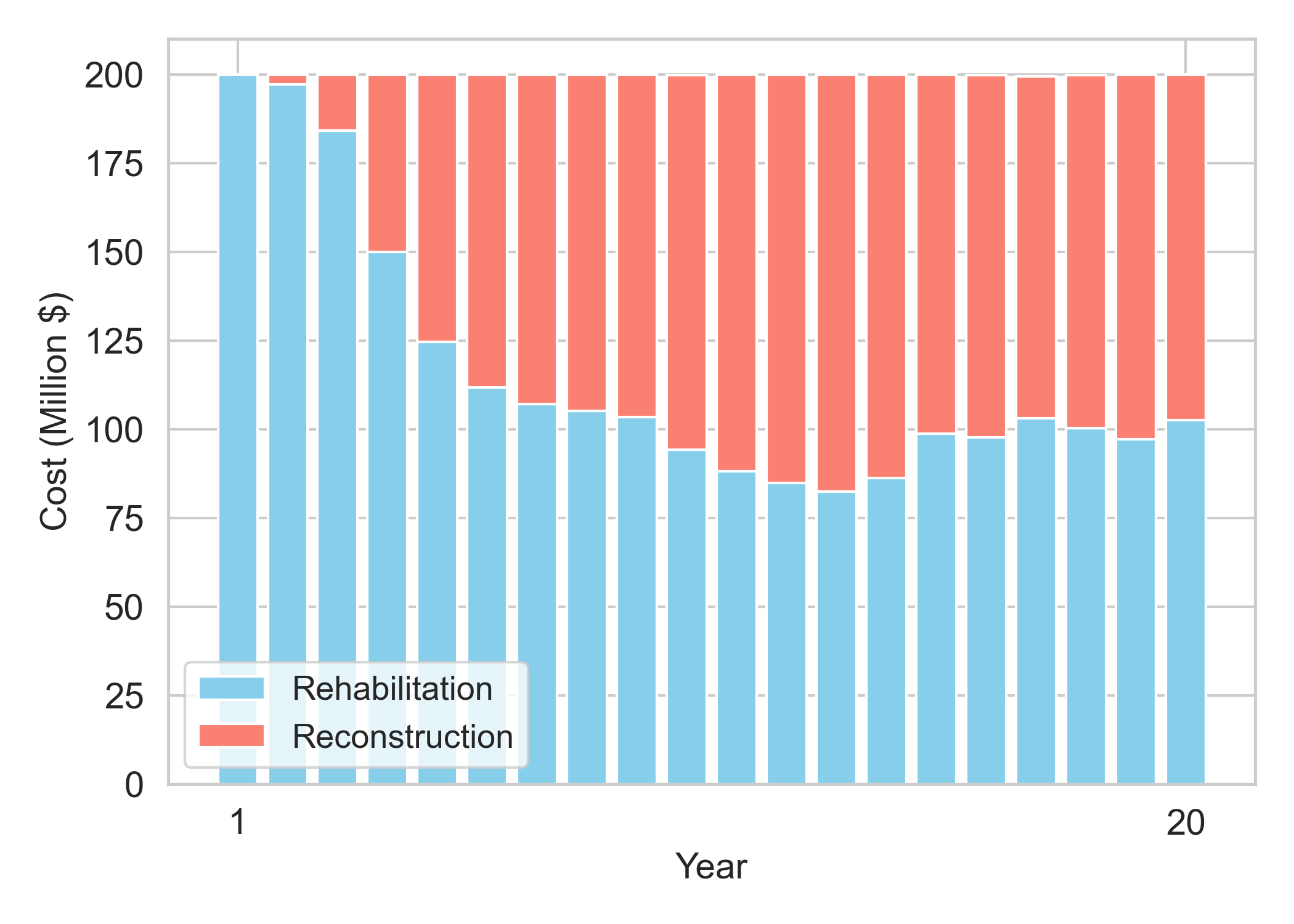}
        \caption{Hybrid LP-GA}
        \label{fig:lp_ga_annual_cost}
    \end{subfigure}
    \hfill
    \begin{subfigure}{0.48\textwidth}
        \includegraphics[width=\textwidth]{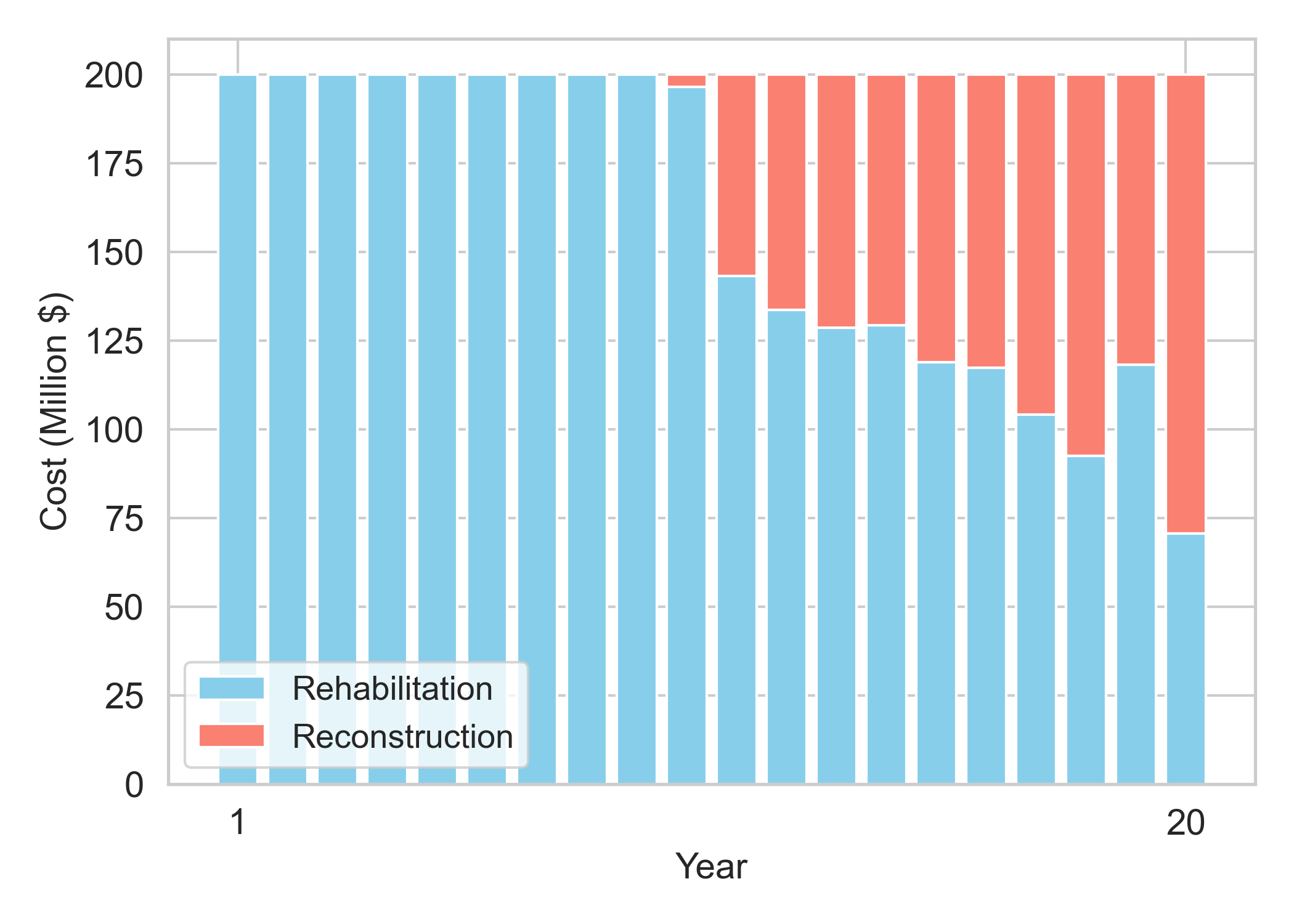}
        \caption{Network DQL}
        \label{fig:drl_annual_cost}
    \end{subfigure}
    
    \caption{Annual cost breakdown for each strategy, illustrating the share of rehabilitation vs. reconstruction.}
    \label{fig:comparison_cost}
\end{figure}
Figure~\ref{fig:comparison_cost} compares the annual maintenance expenditures across the four solution methods, illustrating how each balances rehabilitation versus reconstruction actions over the 20-year horizon. 

Under the Worst-First strategy (Figure~\ref{fig:worst_first_annual_cost}), nearly the entire budget is spent on reconstruction every year. Because high-priority segments are persistently in poor condition, this approach rapidly depletes funds on reconstruction treatments, leaving no fund for assets that need rehabilitation.

In contrast, the Progressive LP method (Figure~\ref{fig:prog_annual_cost}) starts by allocating most of its budget to rehabilitation during the first few years. By around year 3, however, it includes substantial reconstruction expenditures and, from that point onward, reverts roughly half of the annual budget to reconstruction.

The Hybrid LP–GA method (Figure~\ref{fig:lp_ga_annual_cost}) budget allocation pattern is almost similar Progressive LP.  Early in the horizon, the majority of the budget is devoted to rehabilitation, gradually transitioning to more reconstruction and after year 7 nearly half of the budget is allocated to reconstruction interventions.

By comparison, the Network DQL strategy (Figure~\ref{fig:drl_annual_cost}) exhibits a marked preference for rehabilitation in the first 10 years. Only in later stages does it shift portions of the budget to reconstruction, reflecting a learned policy that prioritizes cost-effective improvements early on, deferring costlier interventions until the network condition requires them.

\subsubsection{Delayed Reward}

\begin{figure}[ht]
    \centering
    \includegraphics[width=0.85\textwidth]{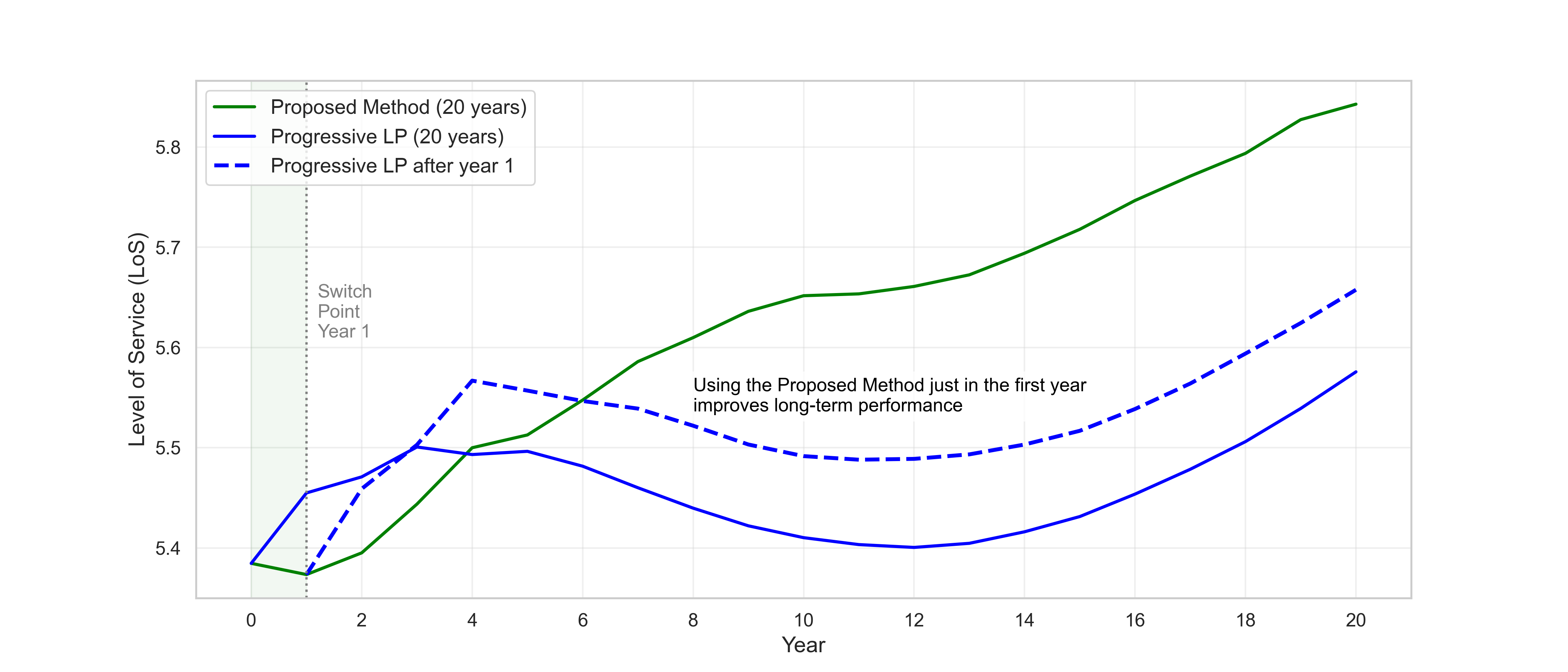}
    \caption{LoS comparison for a hybrid approach (Network DQL in Year~1, followed by Progressive LP in Years~2--20).}
    \label{fig:comparison_hybrid}
\end{figure}

One might have noticed from Figure~\ref{fig:comparison} that the solution of the proposed method looks inferior to the progressive LP and the hybrid LP-GA methods for the first 4 years, although after that the proposed method outperforms the others. To investigate the initial years, we experimented with a hybrid strategy by which the first year's plan takes from the result of the proposed method, and afterwards, the plan is switched back to the progressive LP for the reaming 19 years.  Figure~\ref{fig:comparison_hybrid} compares this hybrid strategy  (dashed blue line) with the Proposed Method (green line) for all 20 years and   Progressive LP (solid blue line) for all 20 years. 

The chart shows that even using the Proposed Method solely in the first year, followed by reversion to Progressive LP, confers a noticeable long-term benefit (the dashed blue curve), compared with Progressive LP alone from the outset (the solid blue curve). In other words, a single year of globally optimized maintenance actions in year~1 helps the network reach a stronger baseline. Although the full-horizon Proposed Method still outperforms this hybrid version, the figure illustrates the value of incorporating globally strategic planning, even if only at the first year, to boost long-term outcomes.

\begin{figure}[ht]
    \centering
    \includegraphics[width=0.85\textwidth]{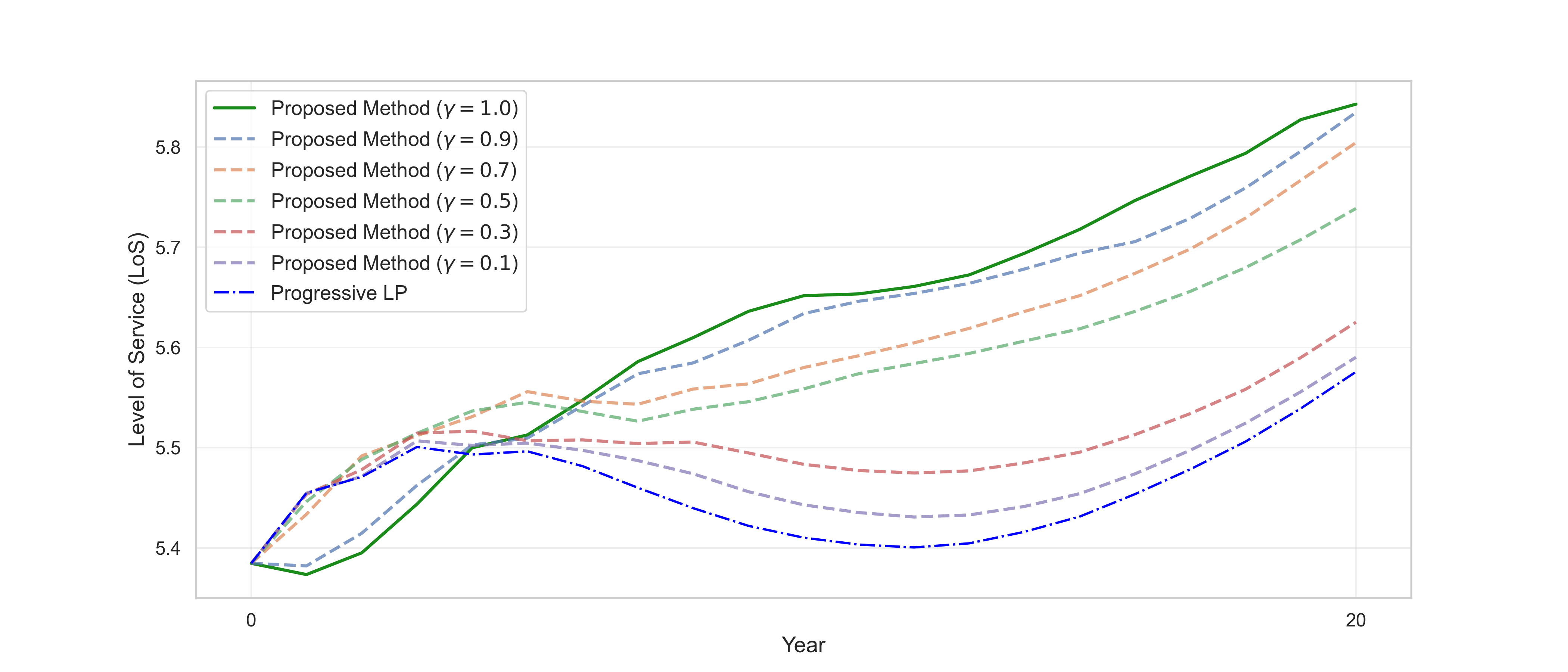}
    \caption{Comparison of different discount factors (\(\gamma\)) on the Proposed Method’s LoS profile over a 20-year horizon, benchmarked against Progressive~LP.}
    \label{fig:comparison_gamma}
\end{figure}

\subsubsection{Effects of Discount Rate}
To investigate the sensitivity of the proposed approach to the discount factor \(\gamma\), Figure~\ref{fig:comparison_gamma} plots the average Level of Service (LoS) trajectories under several choices of \(\gamma\in\{1.0,\,0.9,\,0.7,\,0.5,\,0.3,\,0.1\}\), along with the performance of a Progressive LP benchmark. Here, \(\gamma=1.0\) indicates no discounting of future rewards, giving equal weight to LoS across the entire 20-year planning horizon, whereas lower values of \(\gamma\) bias the decision-making toward more immediate gains in LoS . 

The figure reveals two main observations. First, all variants of the proposed method eventually outperform the Progressive~LP strategy, highlighting the advantage of learning a long-term policy and explicitly planning for future budget needs rather than adopting a purely myopic year-by-year LP. Second, the discount factor substantially influences the timing and magnitude of the LoS gains. Larger \(\gamma\) values (\(\gamma=0.9\) or \(1.0\)) lead to higher EHLoS but show slower initial growth, indicating that the agent defers certain interventions in favor of sustaining long-term performance. By contrast, smaller \(\gamma\) (e.g., \(\gamma=0.3\) or \(0.1\)) produces more immediate improvements but tapers off later in the horizon, as the agent’s policy gives less weight to distant-year outcomes.

From a practical standpoint, infrastructure managers may adopt different discount factors depending on whether their strategic focus is on near-term condition improvements or on sustaining asset performance over many decades. The results in Figure~\ref{fig:comparison_gamma} demonstrate that by tuning \(\gamma\), agencies can fine-tune the trade-off between immediate, incremental benefits and robust, long-term planning.

Overall, these results confirm that the Network DQL method, by considering the long-term impacts of maintenance decisions and tailoring cost allocations over the entire planning period, achieves the highest LoS metrics and offers a viable, scalable solution for large-scale infrastructure management.
\section{Conclusion}
\label{sec:conclusion}

This paper has presented a novel Deep Reinforcement Learning (DRL) framework, referred to as Network DQL, for optimizing maintenance plans of  large-scale infrastructure networks over a multi-year planning horizon. The method tackles two core challenges frequently encountered in large-scale optimization: achieving scalability in the presence of thousands of assets with diverse deterioration behaviors, and adhering to budget constraints that induce intricate interdependencies among maintenance actions. By decomposing the network-level Markov Decision Process (MDP) into asset-level MDPs and incorporating a dedicated budget-allocation mechanism, the approach avoids the exponential growth in computational burden while maintaining annual budget feasibility.

A unified neural network architecture, comprising a Local Q Network, a Policy Network, and a Global Value Network, supports both local and global objectives. The shared-parameter design allows efficient learning across large networks, and the carefully structured state representation, which merges asset-specific features, cost parameters, current conditions, and a compressed view of the entire network, ensures that decisions account for individual asset needs while simultaneously promoting system-level performance goals. Empirical evaluation on a real-world pavement case study with 68{,}800 segments indicates that the proposed method surpasses conventional benchmarks, such as Worst-First strategies, year-by-year linear programming, and a hybrid LP–GA approach, in terms of both Horizon-Averaged LoS and End-of-Horizon LoS.

A principal insight from the experiments is that myopic approaches generally underperform because they fail to anticipate future trade-offs in budget allocation and system-wide condition evolution.

While the proposed framework shows significant promise, several avenues remain open for further research. Enriching the cost model to account for user costs, social impacts, or reliability requirements could offer a more comprehensive assessment of maintenance decisions. Capturing correlated or spatial dependencies, such as spatially correlated deterioration may provide more accurate representations of infrastructure networks. Finally, extending the methodology to robust optimization under deep uncertainty or multi-objective settings would enable the concurrent consideration of long-term life-cycle costs, environmental objectives, and equity among other critical stakeholder priorities.

Overall, the Network DQL approach provides a scalable, data-driven solution for sustainable, long-horizon infrastructure asset management. Its success highlights the potential for DRL-based techniques to address real-world, resource-constrained problems and establishes a framework for broader advances in the strategic planning and upkeep of essential infrastructure systems.

\section{Data Availability}
Some or all data, models, or codes that support the findings of this study are available from the corresponding author upon reasonable request.

\section{Acknowledgements}
The authors gratefully acknowledge the financial support provided by the Natural Sciences and Engineering Research Council of Canada (NSERC) Discovery Grant (RGPIN-2022-04591).

\bibliographystyle{apalike}

\end{document}